\documentstyle[amscd,amssymb,verbatim,diagrams]{amsart}
\newcommand{\w}{{\boldsymbol w}} 

\newcommand{\tr}{\operatorname{tr}}

\newcommand{\Id}{\operatorname{Id}}

\newcommand{\ch}{\operatorname{ch}}

\renewcommand{\mod}{\operatorname{mod}}

\newcommand{\Tr}{\operatorname{Tr}}

\newcommand{\Tor}{\operatorname{Tor}}

\newcommand{\OO}{{\cal O}}

\newcommand{\coker}{\operatorname{coker}}
\newcommand{\DD}{{\cal D}}

\newcommand{\tL}{{\Bbb L}}

\newcommand{\mg}{{\mathfrak m}}

\newcommand{\lan}{\langle}
\newcommand{\ran}{\rangle}

\newcommand{\CC}{{\cal C}}

\newcommand{\si}{\sigma}

\newcommand{\ga}{\gamma}
\newcommand{\de}{\delta}

\numberwithin{equation}{section}

\newtheorem{thm}{Theorem}[section]
\newtheorem{prop}[thm]{Proposition}
\newtheorem{lem}[thm]{Lemma}
\newtheorem{cor}[thm]{Corollary}

{  \theoremstyle{definition}
           \newtheorem{defi}[thm]{Definition}
           \newtheorem{rem}[thm]{Remark}
           \newtheorem{rems}[thm]{Remarks}
          \newtheorem{ex}[thm]{Example}
          
}

\newcommand{\Pf}{\noindent {\it Proof}}
\newcommand{\id}{\operatorname{id}}

\newcommand{\ov}{\overline}
\newcommand{\we}{\wedge}

\newcommand{\rk}{\operatorname{rk}}

\newcommand{\FF}{{\cal F}}
\newcommand{\EE}{{\cal E}}

\newcommand{\SS}{{\cal S}}

\newcommand{\Om}{\Omega}

\newcommand{\Hom}{\operatorname{Hom}}

\newcommand{\Ext}{\operatorname{Ext}}
\newcommand{\End}{\operatorname{End}}

\renewcommand{\a}{\alpha}
\renewcommand{\b}{\beta}

\newcommand{\De}{\Delta}

\newcommand{\C}{{\Bbb C}}

\newcommand{\Z}{{\Bbb Z}}

\newcommand{\Nm}{\operatorname{Nm}}

\newcommand{\Ga}{\Gamma}

\newcommand{\wt}{\widetilde}
\newcommand{\ot}{\otimes}

\newcommand{\sub}{\subset}
\newcommand{\ed}{\qed\vspace{3mm}}

\newcommand{\per}{\operatorname{Per}}
\newcommand{\sTr}{\operatorname{str}}
\newcommand{\sdim}{\operatorname{sdim}}
\newcommand{\MF}{\operatorname{MF}}
\newcommand{\pa}{\partial}
\newcommand{\unit}{{\bf 1}}

\newcommand{\Homb}{{\cal H}om} 

\title{Lefschetz type formulas for dg-categories}
\author{Alexander Polishchuk}
\address{Department of Mathematics, University of Oregon, Eugene, OR 97405}
\email{apolish@@math.uoregon.edu}
\thanks{Supported in part by NSF grant}

\begin{document}
\begin{abstract} We prove an analog of the holomorphic Lefschetz formula for
endofunctors of smooth compact dg-categories. We deduce from it a generalization of
the Lefschetz formula of V.~Lunts \cite{Lunts} that takes the form of a reciprocity law for a pair of commuting
endofunctors. As an application, we prove a version of Lefschetz formula proposed by Frenkel and Ng\^{o}
in \cite{FN}. Also, we compute explicitly the ingredients of the holomorphic Lefschetz formula for the dg-category
of matrix factorizations of an isolated singularity $\w$. We apply this formula to get some restrictions on the Betti
numbers of a $\Z/2$-equivariant module over $k[[x_1,\ldots,x_n]]/(\w)$ in the case when $\w(-x)=\w(x)$.
\end{abstract}
\maketitle

The derived categories of coherent sheaves and their dg-enhancements have been crucial in some recent developments in algebraic geometry. Among these is the new point of view on the circle of ideas related to Chern characters and Hirzebruch-Riemann-Roch theorem, started in the work of Markarian 
\cite{Mark} and continued in \cite{Mark2}, \cite{Cal}, \cite{CW}, \cite{Ram}, \cite{Shk}, \cite{PV-HRR}.
Namely, as shown by Shklyarov in \cite{Shk}, one can formulate and prove
a version of the HRR formula for any smooth and proper dg-category $\CC$,
where the cohomology of the variety is replaced by the Hochschild homology $HH_*(\CC)$. 
The difficulty of applying Shklyarov's formula is usually in the explicit calculation of the Chern character with values in $HH_*(\CC)$ and of the canonical pairing on $HH_*(\CC)$. For example, for derived categories of coherent sheaves on smooth projective varieties this was done in \cite{Cal} and
\cite{Ram}. The case of the category of matrix factorizations associated with an isolated hypersurface singularity was treated in \cite{PV-HRR}.
On the other hand, Lunts proved in \cite{Lunts} a version of the topological Lefschetz trace formula for endofunctors of dg-categories that reduces to the standard one in the case of the derived category
of coherent sheaves. In the present paper we will consider a dg-version of the holomorphic
Lefschetz formula (in its algebraic version) and will calculate its ingredients explicitly for
the category of matrix factorizations of an isolated hypersurface singularity.

Recall that the classical Lefschetz fixed point formula equates the number of fixed points of an endomorphism
$f:M\to M$ of a compact oriented manifold with the supertrace of the action of $f$ on the cohomology of $M$
(assuming the intersection of the graph of $f$ with the diagonal is transversal).
In the case when $M$ is a compact complex manifold and $f$ is holomorphic, the holomorphic Lefschetz formula computes the supertrace of the action of $f$ on $H^*(M,V)$,
where $V$ is a holomorphic vector bundle equipped with a map $f^*V\to V$ (see
\cite[Thm.\ 2]{AB} and \cite{TT}). There are also purely algebraic versions of this theorem (see 
\cite{Donovan}, \cite{Beauv}, \cite{Nielsen}, \cite[II.6]{SGA5}).

In the abstract dg-context the role of the pair $(M,f)$ is played by a pair 
$(\CC,F)$, where 
$\CC$ is a dg-category over a field $k$ and $F$ is an endofunctor of $\per(\CC)$, 
the perfect derived category of $\CC$ (see \cite{Keller-dg}). The cohomology of $M$ is replaced in this
context by the Hochschild homology of $\CC$. The classical situation is recovered when
$\CC$ is a dg-version of the bounded
derived category $D^b(M)$ of coherent sheaves on a smooth projective variety $M$ over
$\C$, and $F$ is the pull-back functor with respect to an endomorphism of $M$. Note that in this
case, by the Kostant-Hochschild-Rosenberg theorem,
the Hochschild homology $HH_*(\CC)$ is isomorphic to the usual cohomology $H^*(M,\C)$.
The assumption of $M$ being smooth and compact has a well-known
analog for dg-categories. In addition we assume the existence of a generator and consider only
endofunctors induced by a kernel in $\per(\CC^{op}\ot\CC)$ (i.e., $F$ is a {\it tensor functor} in the terminology
of \cite{Keller-dg}).

To extend the setup of the holomorphic Lefschetz formula to an arbitrary dg-category 
we interpret $H^*(M,V)$ as $\Ext^*(\OO_M,V)$ and note that $(\OO_M,V)$ can be replaced by any pair
of finite complexes of holomorphic vector bundles (interacting with $f$ appropriately). 
Thus, our abstract formula 
requires as an input a pair of objects $A,B\in\per(\CC)$ equipped with morphisms $\a:A\to F(A)$ and 
$\b:F(B)\to B$. There is an induced endomorphism
\begin{equation}\label{F-a-b-eq}
(F,\a,\b)_*:\Hom(A,B)\rTo{F} \Hom(F(A),F(B))\to\Hom(A,B),
\end{equation}
where the second arrow is given by pre-composing with $\a$ and post-composing with $\b$,
and our formula reads as 
\begin{equation}\label{HLF-eq}
\sTr\bigl((F,\a,\b)_*,\Hom(A,B)\bigr)=\lan\tau^A_F(\a),\tau^{B}_{G}(\wt{\b})\ran_{F,G},
\end{equation}
where 
$G:\per(\CC)\to\per(\CC)$ is the right adjoint functor to $F$,
$\wt{\b}:B\to G(B)$ is induced by $\b$,
$$\tau^A_F:\Hom(A,F(A))\to \Tr_{\CC}(F)$$
is a certain {\it generalized boundary bulk map} taking values in a graded vector space $\Tr_\CC(F)$
associated with $F$, and $\lan ?,?\ran_{F,G}$ is a canonical perfect pairing between
$\Tr_\CC(F)$ and $\Tr_\CC(G)$ that we will define.
More precisely, here $\Tr_\CC$ is the functor on $\CC-\CC$-bimodules
given by tensor multiplication with the diagonal bimodule, so $\Tr_\CC(F)$ is just the Hochschild homology
with values in the bimodule corresponding to $F$. 
In the case when $F=\Id_{\CC}$ the space $\Tr_\CC(F)$ becomes $HH_*(\CC)$ and the formula \eqref{HLF-eq} 
is exactly the generalized abstract Hirzebruch-Riemann-Roch Theorem \cite[Thm.\ 3.1]{PV-HRR} (a similar property is called ``Baggy Cardy Condition" in \cite{CW}).

The data consisting of the spaces $\Tr_\CC(F)$ and $\Tr_\CC(G)$ and the pairing between them
should be thought of as an abstract replacement of the derived intersection of the
graph of an endomorphism $f:M\to M$ with the diagonal. 
In the geometric situation (assuming that the intersection is transversal) the 
spaces $\Tr_\CC(F)$ and $\Tr_\CC(G)$ have bases numbered by the fixed points of $f$, and the
pairing $\lan?,?\ran_{F,G}$ is diagonal with respect to these bases. It is not clear whether one
can formulate an analog of transversality in the abstract situation.

The dg-version of the topological Lefschetz formula, recently established by Lunts \cite{Lunts},
has as an input a pair $(\CC,F)$ as above and states the equality
\begin{equation}\label{LF-eq}
\sTr(F_*,HH_*(\CC))=\sdim\Tr_\CC(F)
\end{equation}
in the ground field $k$,
where $F_*$ is the map on Hochschild homology induced by $F$ and
in the right-hand side we use the superdimension of a graded vector space defined by
$\sdim(V)=\dim(V^{ev})-\dim(V^{odd})$.

Using formula \eqref{HLF-eq} we establish a certain generalization of \eqref{LF-eq} that we call
{\it Lefschetz reciprocity}. Namely, assume that we are given a pair of tensor endofunctors
of $\per(\CC)$, $F$ and $\Psi$, together with a morphism
$$f: F\circ\Psi\to \Psi\circ F.$$
Let $G$ be the right adjoint to $F$. Then
$f$ induces a morphism $\psi:\Psi\circ G\to G\circ\Psi$.
Also, we have the induced map
$$(F,f)_*:\Tr_{\CC}(\Psi)\to \Tr_\CC(G\circ F\circ\Psi)\rTo{G\circ f\circ F} \Tr_\CC(G\circ\Psi\circ F)\rTo{\sim}
\Tr_\CC(F\circ G\circ \Psi)\to \Tr_\CC(\Psi),$$
where 
the third arrow uses the key property that one can switch the order of composition under $\Tr$
(see \eqref{Tr-switch}).
Similarly, we have the map $(\Psi,\psi)_*:\Tr(G)\to\Tr(G)$.
We prove the following reciprocity relation:
\begin{equation}\label{LR-eq0}
\sTr((F,f)_*,\Tr(\Psi))=\sTr((\Psi,\psi)_*,\Tr(G)).
\end{equation}
The idea of the proof is to apply \eqref{HLF-eq}
to a certain endofunctor of $\per(\CC^{op}\ot\CC)$ induced by $F$. 
Note that \eqref{LF-eq} is a particular case of \eqref{LR-eq0}: one should take $\Psi$ to be the
identity functor.


In the case when $F$ is an autoequivalence, we deduce from \eqref{LR-eq0}
a Lefschetz type formula for the action of
an autoequivalence $F$ of $\per(\CC)$ on the Hochschild {\it cohomology} of $\CC$.
Namely, we show that its supertrace is equal to the supertrace of the automorphism
of the space $\Tr(F)$ induced by the inverse Serre functor (see Section \ref{H-coh-sec}).

There are two important classes of examples in which the abstract formulas 
\eqref{HLF-eq} and \eqref{LR-eq0} can be made explicit:
the categories of coherent sheaves on smooth projective varieties and the categories of matrix factorizations of a potential with an isolated singularity.
Note that in \cite{Lunts} the formula \eqref{LF-eq} is rewritten in explicit terms for the case of coherent sheaves
(see also \cite{CT}),
while in \cite{PV-HRR} the case $F=\Id_{\CC}$ of the formula \eqref{HLF-eq} is made explicit
for matrix factorizations. In Section \ref{mf-sec} we extend this calculation to the case of an autoequivalence $F$ induced by a diagonal linear map preserving the potential. In the simplest case
when the origin is an isolated fixed point of this map we get the following result.

\begin{thm}\label{mf-isolated-thm} Let $t=(t_1,\ldots,t_n)\in (k^*)^n$ 
be a diagonal symmetry of an isolated singularity $\w(x_1,\ldots,x_n)\in k[[x_1,\ldots,x_n]]$, such that
$t_i\neq 1$ for all $i$ (i.e., the origin is an isolated fixed point of $t$).
Then for any pair of matrix factorizations of $\w$, $\bar{A}$ and $\bar{B}$, equipped with closed 
morphisms $\a:\bar{A}\to t^*\bar{A}$ and $\b:t^*\bar{B}\to\bar{B}$, one has
$$\sTr((t^*,\a,\b)_*,\Hom(\bar{A},\bar{B}))=\sTr(\a|_0,\bar{A}|_0)\cdot\sTr(\b|_0,\bar{B}|_0)
\cdot\prod_{i=1}^n (1-t_i)^{-1}.$$
\end{thm}

In the case $t=-1$ this formula leads to the divisibility
$$2^{\left \lceil \frac{n}{2} \right  \rceil}|\sTr((-1)^*,\bar{A}|_0)$$
for any $\Z/2$-equivariant matrix factorization $\bar{A}$ of an even isolated singularity $\w$
(see Proposition \ref{divisibility-cor} and Example \ref{divisibility-ex}). In the case of
the matrix factorization obtained by stabilization of a free resolution of a graded module $M$ over 
$k[x_1,\ldots,x_n]/(\w)$ we get restrictions on the Betti numbers of $M$ which are similar but different
from those obtained in \cite{AB-bounds}.

On the other hand, we apply the formula 
\eqref{LR-eq0} in the context of the derived category of coherent sheaves 
on a smooth projective variety
to prove a version of holomorphic Lefschetz formula 
conjectured by Frenkel and Ng\^{o} in
\cite[Conj.\ 6.2]{FN} (see Theorem \ref{der-int-thm}).
Presumably a similar proof should work in a more general
setting of derived stacks, as envisioned in \cite{FN}.

\medskip

\noindent
{\it Convention}. In various maps involving tensor products the Koszul sign rule is tacitly assumed.
We work over a fixed ground field $k$. Whenever we talk about varieties over $k$ we assume that $k$
is algebraically closed. For a morphism of schemes $f$ we denote by $f^*$ and $f_*$ the corresponding
derived functors.

\noindent
{\it Acknowledgments}. I am grateful to Luchezar Avramov and David Eisenbud for helpful discussions.


\section{Kernels and traces}\label{trace-sec}

\subsection{Kernels, duals and tensor products}

We refer to \cite{Keller-derived}, \cite{Keller-dg}, \cite{Toen-dg} and \cite{TV}
for results about dg-modules over dg-categories.
For a dg-category $\CC$ we denote by $\per(\CC)$ the derived category of perfect $\CC^{op}$-modules, and by 
$\per_{dg}(\CC)$ its dg-version, so that $\per(\CC)$ is the homotopy category of $\per_{dg}(\CC)$.
For $M,N\in\per(\CC)$ we denote by $\Hom(M,N)=\Hom_{\CC^{op}}(M,N)$ the graded
morphism spaces in $\per(\CC)$.
Note that these are given by the cohomology of the complexes
$\Hom_{\CC^{op}-\mod}(M,N)$.
In this paper we consider only dg-categories $\CC$ such that 
$\per_{dg}(\CC)$ is saturated, i.e., smooth, compact and has a generator.
This implies that the spaces $\Hom(M,N)$ are finite-dimensional.

Every dg-functor $F:\per_{dg}(\CC)\to\per_{dg}(\DD)$ (up to an appropriate equivalence)
can be given by a kernel $K\in\per(\CC^{op}\ot\DD)$, so that
$$F=T_K: M\mapsto M\ot_{\CC}^{\tL} K.$$
This correspondence can be formulated as a certain quasi-equivalence of dg-categories 
(see \cite[Sec.\ 2.2]{TV}, \cite[Thm.\ 4.6]{Keller-dg}). 
In this paper, whenever we talk about a functor $F:\per(\CC)\to\per(\DD)$, we implicitly assume that it
is given on a dg-level, and hence is induced by a kernel. Sometimes, we denote the corresponding
kernel in $\per(\CC^{op}\ot\DD)$ also by $F$.
We denote by $\De_{\CC}\in\per(\CC^{op}\ot\CC)$ the diagonal kernel representing the identity functor
on $\per(\CC)$.

The composition of dg-functors $F_1\circ F_2$, where $F_2$ is given by a kernel $K_2\in\per(\CC^{op}\ot\DD)$
and $F_1$ is given by a kernel $K_1\in\per(\DD^{op}\ot\EE)$,
corresponds to the kernel $K_2\ot_{\DD}^{\tL} K_1\in\per(\CC^{op}\ot\EE)$. 
Sometimes, when $K_i$ are implicit, we denote this kernel simply by $F_1\circ F_2$. 

Recall that we have a natural duality
\begin{equation}\label{duality-functor}
\per(\CC)^{op}\rTo{\sim}\per(\CC^{op}):M\mapsto M^\vee,
\end{equation}
where $M^\vee(C)=\Hom_{\CC^{op}-\mod}(M,h_C)$ (morphisms in the dg-category of $\CC^{op}$-modules). Here $h_C$ is the representable $\CC^{op}$-module associated with $C\in\CC$. 
For $M\in\per_{dg}(\CC)$ and $N\in\per_{dg}(\CC^{op})$ we have a natural quasi-isomorphism 
\begin{equation}\label{M-N-dual-eq}
M\ot_{\CC}^{\tL} N^\vee\to\Hom_{\CC^{op}-\mod}(N,M)
\end{equation}
(see \cite[6.2]{Keller-derived}). Also, we have
a natural quasi-isomorphism $M\to (M^\vee)^\vee$.

For kernels
$K_1\in\per(\CC_1^{op}\ot\DD_1)$ and $K_2\in\per(\CC_2^{op}\ot\DD_2)$
we set
$$K_1\Box K_2=(K_1\ot K_2)\circ \si_{23}\in\per(\CC_1^{op}\ot \CC_2^{op}\ot\DD_1\ot\DD_2),$$
where $\si_{23}$ is the permutation.
We denote the corresponding functor by
$$F_1\Box F_2:\per(\CC_1\ot\CC_2)\to\per(\DD_1\ot\DD_2),$$
where $F_1$ and $F_2$ are functors associated with $K_1$ and $K_2$.

\subsection{Adjoint kernels}\label{adj-sec}

Let 
$$F=T_K:\per(\CC)\to\per(\DD): M\mapsto M\ot_\CC K$$ be 
the functor associated with a kernel $K\in\per(\CC^{op}\ot\DD)$.
We associate with $K$ two other kernels. First, we define
$K^T\in\per(\DD^{op}\ot\CC)$, the {\it right dual kernel to} $K$, by
$$K^T(D,C^\vee)=\Hom_{\DD^{op}-\mod}(K(C,?),h_D),$$
where $h_D$ is the representable right $\DD$-module associated with $D$.
The corresponding functor 
$$G=T_{K^T}:\per(\DD)\to\per(\CC)$$ 
is right adjoint to $F$
(see \cite[Sec.\ 6.2]{Keller-derived} and \cite[Sec.\ 1.2]{PV-HRR}).
Furthemore, the adjunction morphisms are induced by the canonical morphisms of kernels
in $\per(\DD^{op}\ot\DD)$ and $\per(\CC^{op}\ot\CC)$:
\begin{equation}\label{adj-to-id}
K^T\ot_\CC^{\tL} K\to\De_{\DD} \ \text{ and}
\end{equation}
\begin{equation}\label{adj-from-id}
\De_{\CC}\to K\ot_{\DD}^{\tL} K^T
\end{equation}
(in fact, the first morphism has a natural realization in $\per_{dg}(\DD^{op}\ot\DD)$, see
\cite[(1.14)]{PV-HRR}).
In terms of the duality \eqref{duality-functor} we can rewrite the definition of $K^T$ as follows:
$$K^T(?,C^\vee)=K(C,?)^\vee$$

We also set 
$$K':=K\circ\si\in\per(\DD\ot\CC^{op}),$$ 
where $\si$ is the permutation of factors, and denote by
$F':\per(\DD^{op})\to\per(\CC^{op})$ the corresponding functor.

\begin{lem}\label{G'-F'-lem} 
The pair $(G',F')$ is also adjoint. More precisely, there is a natural quasi-isomorphism
of kernels $K'\to (K^T)^{\prime T}$.
\end{lem}

\Pf . The functor $G'$ is given by
the kernel $(K^T)'=K^T\circ\si\in\per(\CC\ot\DD^{op})$ satisfying
$$(K^T)'(C^\vee,?)=K(C,?)^\vee.$$
Hence, the right adjoint to $G'$ is given by the kernel 
$(K^T)^{\prime T}\in\per(\DD\ot\CC^{op})$. 
We have
$$(K^T)^{\prime T}(?,C)=(K^T)'(C^\vee,?)^\vee=K(C,?)^{\vee\vee}$$
and there is a natural quasi-isomorphism $K'(?,C)=K(C,?)\to K(C,?)^{\vee\vee}$.
\ed

\begin{lem}\label{F-op-lem} 
(i) For $M\in\per(\DD)$ and $N\in\per(\CC)$ one has a natural isomorphism
$$\Hom_{\DD^{op}}(M,F(N))\simeq\Hom_{\CC}(N^\vee,F'(M^\vee)).$$

\noindent
(ii) For $N\in\per(\CC)$ there is a natural isomorphism
\begin{equation}\label{G'-F-dual-eq}
G'(N^\vee)\simeq F(N)^{\vee}
\end{equation}
in $\per(\DD^{op})$, so that the following diagram of isomorphisms is commutative
\begin{equation}\label{dual-adjoint-diagram}
\begin{diagram}
\Hom_{\DD^{op}}(M,F(N))&\rTo{}&\Hom_{\CC}(N^\vee,F'(M^\vee))\\
\dTo{}&&\dTo{\ga}\\
\Hom_{\DD}(F(N)^{\vee},M^{\vee})&\rTo{}&\Hom_{\DD}(G'(N^\vee),M^\vee)
\end{diagram}
\end{equation}
for $M\in\per(\DD)$ and $N\in\per(\CC)$. In this diagram the upper horizontal arrow
is an isomorphism of part (i), the lower horizontal arrow is induced by \eqref{G'-F-dual-eq},
and $\ga$ is the adjunction isomorphism of Lemma \ref{G'-F'-lem}.
\end{lem}

\Pf . (i) By \eqref{M-N-dual-eq}, we have
$$\Hom(M,F(N))\simeq F(N)\ot_{\DD}^{\tL} M^\vee\simeq N\ot_{\CC}^{\tL} K\ot_{\DD}^{\tL} M^\vee.$$
Reversing the order of factors we can rewrite this as
$$M^\vee\ot_{\DD^{op}}^{\tL} K'\ot_{\CC^{op}}^{\tL} N\simeq F'(M^\vee)\ot_{\CC^{op}}^{\tL} N\simeq
\Hom(N^\vee, F'(M^\vee)).$$

\noindent
(ii) First, taking cofibrant resolutions for $K'$ and $(K^T)'$,
we can realize the map $\ga$ (with $M^\vee$ replaced by any $\DD$-module $L$)
on the chain level as the following composition:
\begin{align}\label{ga-comp-eq}
&\Hom(N^\vee,L\ot_{\DD^{op}}K')\to
\Hom(N^\vee\ot_{\CC^{op}} (K^T)', L \ot_{\DD^{op}}K'\ot_{\CC^{op}} (K^T)')\to\nonumber\\
&\Hom(N^\vee\ot_{\CC^{op}} (K^T)', L),
\end{align}
where the second arrow is induced by the map
$$K'\ot_{\CC^{op}}(K^T)'\to (K^T)^{\prime T}\ot_{\CC^{op}} (K^T)'\to\De_{\DD^{op}}$$
in $\per_{dg}(\DD\ot\DD^{op})$. Composing \eqref{ga-comp-eq} with a quasi-isomorphism of the form
$?\ot_{\CC^{op}} N\to\Hom(N^\vee,?)$ obtained from
\eqref{M-N-dual-eq},
we get a quasi-isomorphism
$$L\ot_{\DD^{op}}K'\ot_{\CC^{op}} N\to \Hom(N^\vee\ot_{\CC^{op}} (K^T)', L).$$
Note that $K'\ot_{\CC^{op}} N\simeq N\ot_{\CC} K$. Thus,
for $L=h_{D^\vee}$, where $D^\vee\in\DD^{op}$, this gives a quasi-isomorphism of dg-modules
$$(N\ot_\CC K)(D^\vee)\to (N^\vee\ot_{\CC^{op}} (K^T)')^\vee(D^\vee),$$
hence we get a quasi-isomorphism $F(N)\to G'(N^\vee)^\vee$, from which \eqref{G'-F-dual-eq} is
obtained by duality. 

Now using the isomorphisms \eqref{M-N-dual-eq}
we can replace \eqref{dual-adjoint-diagram} with the following diagram in which we replaced
$M^\vee$ by $L$ and inserted an additional diagonal arrow:
\begin{equation}
\begin{diagram}
N\ot_{\CC} K\ot_{\DD} L&\rTo{}& L\ot_{\DD^{op}} K'\ot_{\CC^{op}} N\\
\dTo{}&\ruTo{}&\dTo{}\\
L\ot_{\DD^{op}}(N\ot_{\CC} K)&\rTo{}&L\ot_{\DD^{op}}G'(N^\vee)^\vee
\end{diagram}
\end{equation}
It remains to observe that in this diagram the upper left triangle is commutative for trivial reasons, and
the lower right triangle is commutative by the above construction of the map $F(N)\to G'(N^\vee)^\vee$.
\ed

Let us consider the external tensor product 
$K'\Box K^T\in\per(\DD\ot \DD^{op}\ot\CC^{op}\ot\CC)$.
Then we have a natural isomorphism
$$(K'\Box K^T)\ot^{\tL}_{\CC^{op}\ot\CC}\De_{\CC^{op}}
\rTo{\sim} (K^T\ot_{\CC} K)\circ\si
$$
in $\per(\DD^{op}\ot\DD)$ (cf. \cite[(1.8)]{PV-HRR}).
Combining it with the permutation of \eqref{adj-to-id} we get a map of kernels
\begin{equation}\label{quasi-adj-map-bis}
(K'\Box K^T)\ot^{\tL}_{\CC^{op}\ot\CC}\De_{\CC^{op}}
\to\De_{\DD^{op}} 
\end{equation}
in $\per(\DD\ot\DD^{op})$. 
Note also that the functor $F'\Box G:\per(\DD^{op}\ot\DD)\to\per(\CC^{op}\ot\CC)$ given
by the kernel $K'\Box K^T$, viewed as an operation on functors, acts by
$\Phi\mapsto G\circ\Phi\circ F$. In particular, $(F'\Box G)(\De_\DD)$ is the kernel giving $G\circ F$.

\subsection{Connection with traces}

Recall that we have the natural trace functors
$$\Tr_{\CC}:\per(\CC^{op}\ot\CC)\to\per(k),$$ 
$$\Tr^{dg}_{\CC}:\per_{dg}(\CC^{op}\ot\CC)\to\per_{dg}(k)$$
given by the tensor product on the right with the diagonal $\CC^{op}\ot\CC$-module $\De_{\CC^{op}}$
(see \cite[Sec.\ 5.2.3]{Toen-dg} and \cite[Sec.\ 1.1]{PV-HRR}). 
In other words, $\Tr_\CC$ is the functor of Hochschild homology with coefficients in a bimodule.
By abuse of notation, for a tensor endofunctor $F:\per(\CC)\to\per(\CC)$ given by a kernel
$K_F\in\per(\CC^{op}\ot\CC)$ we will write $\Tr_\CC(F)$ instead of $\Tr_\CC(K_F)$.
Recall that for any tensor functors $F_1:\per(\CC)\to\per(\DD)$ and $F_2:\per(\DD)\to\per(\CC)$
we have a canonical isomorphism
\begin{equation}\label{Tr-switch}
\Tr_\DD(F_1\circ F_2)\simeq \Tr_\CC(F_2\circ F_1)
\end{equation}
(see \cite[Lem.\ 1.1.3]{PV-HRR}).

Returning to the situation of Section \ref{adj-sec} we observe that the map of kernels \eqref{quasi-adj-map-bis} induces a natural transformation
\begin{equation}\label{trace-nat-transformation}
t(F,G):\Tr_{\CC}\circ (F'\Box G)\to  \Tr_{\DD}.
\end{equation}

By Lemma \ref{G'-F'-lem}, the above construction is also applicable to the adjoint pair of functors
$(G',F')$, so it gives a map
$$t(G',F'):\Tr_{\CC^{op}}\circ (G\Box F')\to\Tr_{\DD^{op}}.$$
It is easy to see that the following diagram of functors
$\per(\DD\ot\DD^{op})\to \per(k)$ is commutative 
\begin{equation}\label{Tr-F-G-comm-diag}
\begin{diagram}
\Tr_{\CC}\circ (F'\Box G)\circ \si_{\DD}&\rTo{t(F,G)\circ\si_{\DD}}&\Tr_{\DD}\circ\si_{\DD}\\
\dTo{\sim}&&\dTo{\sim}\\
\Tr_{\CC^{op}}\circ (G\Box F')&\rTo{t(G',F')}&\Tr_{\DD^{op}}
\end{diagram}
\end{equation}
where $\si_{\DD}:\per(\DD\ot\DD^{op})\to\per(\DD^{op}\ot\DD)$ is the permutation equivalence,
the isomorphism in the left column is induced by the isomorphism
$$(F'\Box G)\circ\si_\DD\simeq \si_{\CC}\circ (G\Box F')$$
and by the isomorphism $\Tr_{\CC}\circ\si_{\CC}\simeq\Tr_{\CC^{op}}$.

Note that we have a natural isomorphism
$$(F'\Box G)(\De_{\DD})\simeq \De_{\DD}\otimes_{\DD^{op}\ot\DD}(K'\Box K^T)\simeq
K\ot_{\DD} K^T.$$
Hence, the canonical adjunction map \eqref{adj-from-id} can be viewed as a map 
\begin{equation}\label{ga-F-G0}
\ga_{F,G}:\De_{\CC}\to (F'\Box G)(\De_{\DD}).
\end{equation}
Using isomorphism of Lemma \ref{F-op-lem}(i) for the map
$$\ga_{G',F'}:\De_{\CC^{op}}\to (G\Box F')(\De_{\DD^{op}})$$
we obtain a canonical morphism
\begin{equation}\label{Serre-morphism0}
\De_{\DD^{op}}^{\vee}\to (G'\Box F)(\De_{\CC^{op}}^{\vee}).
\end{equation}

Note also that for $M\in\per(\CC^{op}\ot\CC)$ we have a canonical isomorphism
\begin{equation}\label{Serre-Tr0}
\Hom(\De_{\CC^{op}}^{\vee},M)\simeq M\ot^{\tL}_{\CC^{op}\ot\CC}\De_{\CC^{op}}=\Tr_{\CC}(M).
\end{equation}

\begin{lem}\label{t-F-G-lem}
(i) For $\Phi\in\per(\DD^{op}\ot\DD)$
the natural transformation 
$$t(F,G)(\Phi):\Tr_{\CC}((F'\Box G)(\Phi))\to\Tr_{\DD}(\Phi)$$
is given by the composition
$$\Tr_{\CC}((F'\Box G)(\Phi))\simeq\Tr_{\CC}(G\circ\Phi\circ F)\simeq\Tr_{\DD}(F\circ G\circ \Phi)\to \Tr(\Phi),$$
where the second isomorphism uses \eqref{Tr-switch}.

\noindent (ii) Via the identification \eqref{Serre-Tr0}, $t(F,G)(\Phi)$ can be identified with
the composition
\begin{equation}\label{t-F-G-Phi-comp}
\Hom(\De_{\CC^{op}}^{\vee},(F'\Box G)(\Phi))\rTo{\sim}
\Hom((G'\Box F)(\De_{\CC^{op}}^{\vee}),\Phi)\to\Hom(\De_{\DD^{op}}^{\vee},\Phi),
\end{equation}
where the first arrow is given by adjunction of the pair $((G'\Box F),(F'\Box G))$, while the second
is induced by \eqref{Serre-morphism0}.
\end{lem}

\Pf . (i) Consider the object 
$$L=K\boxtimes\Phi\boxtimes K^T\in \per(\CC^{op}\ot\DD\ot\DD^{op}\ot\DD\ot\DD^{op}\ot\CC).$$
We have
\begin{equation}\label{Tr-G-Phi-F}
\Tr_{\CC}(G\circ\Phi\circ F)=\Tr_{\CC}(K\ot_{\DD}\Phi\ot_{\DD} K^T)=
\Tr_{\CC}((\De^{23}_{\DD}\boxtimes\De^{45}_{\DD})\ot_{\DD^{(4)}}L),
\end{equation}
where $\DD^{(4)}=\DD^{op}\ot\DD\ot\DD^{op}\ot\DD$.
On the other hand,
\begin{equation}\label{Tr-F-G-Phi}
\Tr_{\DD}(F\circ G\circ\Phi)=\Tr_{\DD}(\De^{45}_{\DD}\ot_{\DD^{op}\ot\DD} L\ot_{\CC^{op}\ot\CC}
\De^{16}_{\CC^{op}}).
\end{equation}
Now the isomorphism between \eqref{Tr-F-G-Phi} and \eqref{Tr-G-Phi-F} is obtained by identifying
both with
$$(\De^{23}_{\DD}\boxtimes\De^{45}_{\DD})\ot_{\DD^{(4)}} L \ot_{\CC^{op}\ot\CC}\De^{16}_{\CC^{op}}.$$
The map $\Tr_{\DD}(F\circ G\circ\Phi)\to\Tr_{\DD}(\Phi)$ is induced by the adjunction morphism
$F\circ G\to\De_{\DD}$, i.e., in terms of the above isomorphism by the map
$$L  \ot_{\CC^{op}\ot\CC}\De^{16}_{\CC^{op}}\to (\Phi\boxtimes\De_{\DD})\circ \si_{4321}$$
induced by \eqref{adj-to-id},
where $\si_{4321}$ is the cyclic permutation of factors in the tensor product.
Now the required compatibility follows from the commutativity of the diagram
\begin{diagram}
\Phi\ot_{\DD^{op}\ot\DD}(K'\Box K^T)\ot_{\CC^{op}\ot\CC}\De_{\CC^{op}}&\rTo{}&
\Phi\ot_{\DD^{op}\ot\DD}\De_{\DD^{op}}\\
\dTo{\sim}&&\dTo{\sim}\\
(\De^{23}_{\DD}\boxtimes\De^{45}_{\DD})\ot_{\DD^{(4)}} L \ot_{\CC^{op}\ot\CC}\De^{16}_{\CC^{op}}&\rTo{}&
(\De^{23}_{\DD}\boxtimes\De^{45}_{\DD})\ot_{\DD^{(4)}}[(\Phi\boxtimes\De_{\DD})\circ \si_{4321}]
\end{diagram}
in which both horizontal arrows are induced by \eqref{adj-to-id}.

\noindent
(ii) Applying the commutative diagram of Lemma \ref{F-op-lem}(ii) to the adjoint pair $((G'\Box F),(F'\Box G))$ and
the pair of objects $(\De_{\CC^{op}}, \De_{\DD^{op}})$ we obtain that the dual map to 
\eqref{Serre-morphism0},
$$(F\Box G')(\De_{\CC^{op}})\simeq \bigl((G'\Box F)(\De_{\CC^{op}}^{\vee})\bigr)^\vee\to\De_{\DD^{op}}$$
corresponds by adjointness to $\ga_{G',F'}$.
Hence, this map is given by the natural morphism
$$(F\Box G')(\De_{\CC^{op}})\simeq (K'\Box K^T)\ot_{\CC^{op}\ot\CC}\De_{\CC^{op}}\to \De_{\DD^{op}}$$
(see \eqref{quasi-adj-map-bis}). Now rewriting the composition \eqref{t-F-G-Phi-comp} using the isomorphisms of the form
\eqref{M-N-dual-eq} we obtain the result.
\ed

\section{Generalized functoriality for Hochschild homology with coefficients}\label{gen-funct-sec}

\subsection{Generalized functoriality}\label{funct-sec}

Let $(\CC_1,F_1)$ and $(\CC_2,F_2)$ 
be a pair of dg-categories equipped with endofunctors $F_i$ of $\per_{dg}(\CC_i)$.

\begin{defi}\label{dg-functor-end-def} 
(i) A {\it dg-functor} 
\begin{equation}\label{Phi-phi-eq}
(\Phi,\phi):(\per_{dg}(\CC_1),F_1)\to (\per_{dg}(\CC_2),F_2)
\end{equation}
consists of a dg-functor $\Phi:\per_{dg}(\CC_1)\to\per_{dg}(\CC_2)$ 
together with a (closed) morphism of dg-functors 
$$\phi:\Phi\circ F_1\to F_2\circ \Phi.$$

\noindent
(ii) If $(\CC_3,F_3)$ is another dg-category with an endofunctor of $\per_{dg}(\CC_3)$ and
$$(\Phi',\phi'):(\per_{dg}(\CC_2),F_2)\to (\per_{dg}(\CC_3),F_3)$$ 
is a dg-functor, then we have an induced morphism
$$\phi'\circ\phi:\Phi'\circ\Phi\circ F_1\to \Phi'\circ F_2\circ\Phi\to F_3\circ \Phi'\circ\Phi$$
and we define the composition of these dg-functors by
$$(\Phi',\phi')\circ (\Phi,\phi)=(\Phi'\circ\Phi,\phi'\circ\phi):(\per_{dg}(\CC_1),F_1)\to 
(\per_{dg}(\CC_3),F_3).$$ 
\end{defi}

Given a dg-functor \eqref{Phi-phi-eq}
we define the induced morphism
\begin{equation}\label{functoriality-map}
(\Phi,\phi)_*:\Tr_{\CC_1}(F_1)\to\Tr_{\CC_2}(F_2)
\end{equation}
as the composition 
$$\Tr(F_1)\to\Tr(\Psi\circ\Phi\circ F_1)\to\Tr(\Phi\circ F_1\circ\Psi)\to\Tr(F_2\circ\Phi\circ\Psi) \to \Tr(F_2),$$
where $\Psi$ is the right adjoint functor to $\Phi$.

The above construction is compatible with compositions. 

\begin{lem}\label{fun-comp-lem} 
In the situation of Definition \eqref{dg-functor-end-def}(ii) one has 
$$(\Phi'\circ\Phi,\phi'\circ\phi)_*=(\Phi',\phi')_*\circ (\Phi,\phi)_*$$
in $\Hom_k(\Tr_{\CC_1}(F_1),\Tr_{\CC_3}(F_3))$.
\end{lem}

\Pf . This follows from the commutativity of the diagram (where we omitted the $\circ$ symbols)
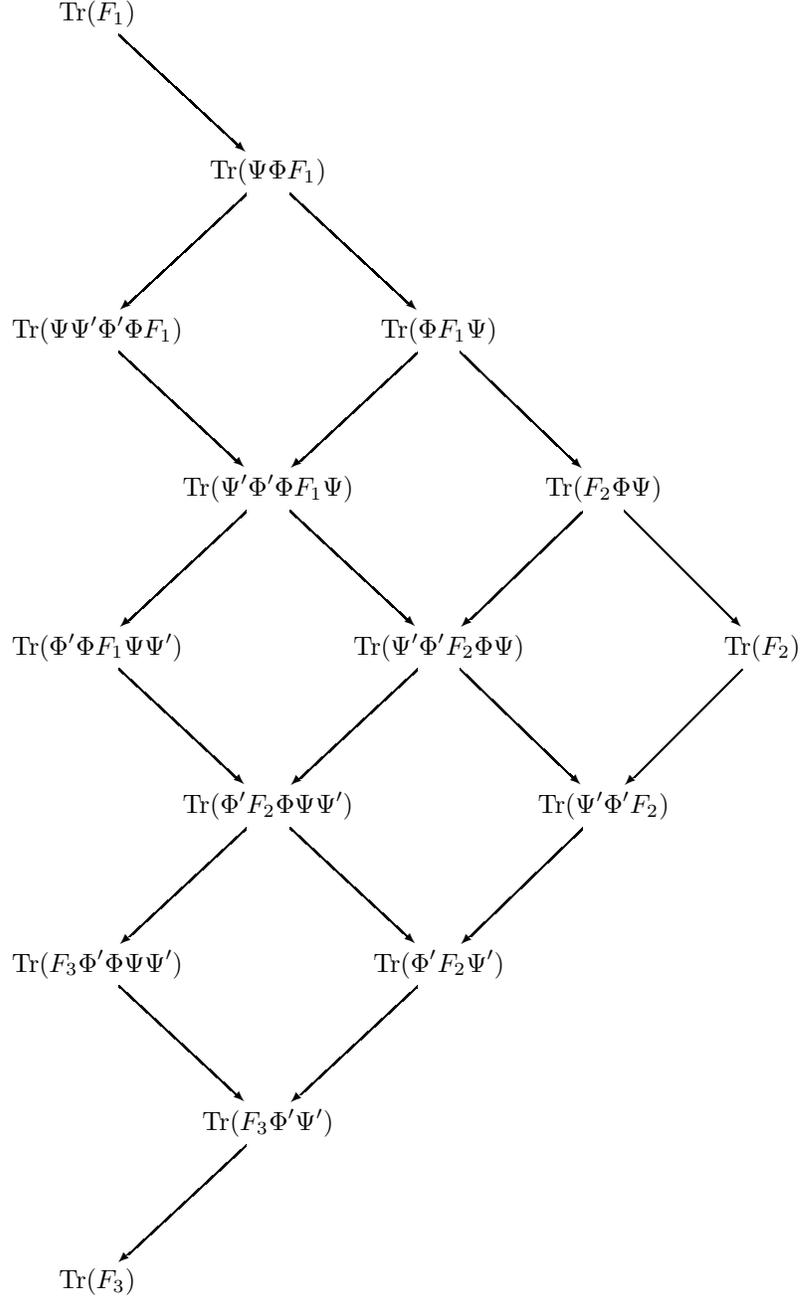
\begin{diagram}
\Tr(F_1)\\
&\rdTo{}\\
&&\Tr(\Psi\Phi F_1)\\
&\ldTo{}&&\rdTo{}\\
\Tr(\Psi\Psi'\Phi'\Phi F_1)&&&&\Tr(\Phi F_1\Psi)\\
&\rdTo{}&&\ldTo{}&&\rdTo{}\\
&&\Tr(\Psi'\Phi'\Phi F_1\Psi)&&&&\Tr(F_2\Phi\Psi)\\
&\ldTo{}&&\rdTo{}&&\ldTo{}&&\rdTo{}\\
\Tr(\Phi'\Phi F_1\Psi\Psi')&&&&\Tr(\Psi'\Phi' F_2\Phi\Psi)&&&&\Tr(F_2)\\
&\rdTo{}&&\ldTo{}&&\rdTo{}&&\ldTo{}\\
&&\Tr(\Phi' F_2\Phi\Psi\Psi')&&&&\Tr(\Psi'\Phi' F_2)\\
&\ldTo{}&&\rdTo{}&&\ldTo{}\\
\Tr(F_3\Phi'\Phi\Psi\Psi')&&&&\Tr(\Phi' F_2\Psi')\\
&\rdTo{}&&\ldTo{}\\
&&\Tr(F_3\Phi'\Psi')\\
&\ldTo{}\\
\Tr(F_3)
\end{diagram}
Namely, $(\Phi'\circ\Phi,\phi'\circ\phi)_*$ is equal to the composition of the arrows going from $\Tr(F_1)$
to $\Tr(F_3)$ along the left edge of the diagram, while $(\Phi',\phi')_*\circ (\Phi,\phi)_*$
is the composition of the arrows along the right edge of the diagram.
\ed

\begin{rems}
1. In the case when $F_1$ and $F_2$ are the identity functors and $\phi=\id$, the
associated map \eqref{functoriality-map} is the functoriality map 
$\Phi_*:HH_*(\CC_1)\to HH_*(\CC_2)$ defined in \cite[Sec.\ 1.2]{PV-HRR}.
The equivalence of this definition of $F_*$ with the standard definition of the functoriality using Hochschild complexes is checked in \cite[Appendix]{PV-CohFT}.

\noindent
2. The definition of the maps $(\Phi,\phi)_*$ is the particular case of the trace maps for a bicategory with
a shadow functor (see \cite[Def.\ 5.1]{PS}). Namely, we have a bicategory whose $0$-cells are dg-categories,
the categories of $1$-cells are given by derived categories of bimodules, and the shadow functors are
the functors $\Tr_{\CC}$ (cf. \cite[Ex.\ 6.4]{PS}). In this context our Lemma \ref{fun-comp-lem} follows from
\cite[Prop.\ 7.5]{PS}.
\end{rems}

\subsection{Generalized boundary-bulk map}\label{bbm-sec}

Given a tensor functor $F:\per(\CC)\to\per(\CC)$ (given by a kernel $K\in\per(\CC^{op}\ot\CC)$)
and an object $A\in\per(\CC)$, we can construct a natural map
$$\tau^A_F:\Hom(A,F(A))\to \Tr_\CC(F)$$
which in the case $F=\Id$ becomes the boundary-bulk map defined in \cite[Prop.\ 1.2.4]{PV-HRR}.
Namely, this map is defined by applying the functor $\Tr_\CC$ to the natural morphism
\begin{equation}\label{A-vee-F-A-eq}
c^A_F: A^\vee\boxtimes F(A)\to F
\end{equation}
in $D(\CC^{op}\ot\CC)$ (see \cite[Eq.\ (1.25)]{PV-HRR}).
It is useful to note that the map $c^A_F$ is obtained by applying $\Id_{\CC}\Box F$ to the morphism
\begin{equation}\label{A-vee-Id-eq}
c^A_{\Id}: A^\vee\boxtimes A\to\De_{\CC},
\end{equation}
or equivalently, by composing the corresponding functors with $F$ on the left
(see the proof of Lemma 1.2.5(ii) in \cite{PV-HRR}).

Note that we can view an object $A$ as a dg-functor 
$\iota_A:\per_{dg}(k)\to \per_{dg}(\CC)$ sending $k$ to $A$.
Then a map $\a:A\to F(A)$ is the same as a morphism of dg-functors $\a:\iota_A\to F\circ\iota_A$.
Thus, we have a dg-functor 
$$(\iota_A,\a):(\per_{dg}(k),\Id)\to (\per_{dg}(\CC),F)$$
and the construction of \ref{funct-sec}
gives an induced map
$$(\iota_A,\a)_*:k\to \Tr_{\CC}(F).$$

\begin{lem}\label{bb-funct-lem} 
One has
$$\tau^A_F(\a)=(\iota_A,\a)_*(1).$$
\end{lem}

\Pf . The right adjoint functor $\Psi$ to $\Phi=\iota_A$ is given by the kernel $A^\vee$, viewed as
a $\CC-k$-bimodule, so $\Psi\circ\Phi$ corresponds to the vector space $\Hom(A,A)$ while
$\Phi\circ\Psi$ corresponds to the kernel $A^\vee\boxtimes A$.
Thus, by definition, $(\iota_A,\a)_*$ is the following composition
$$k\rTo{\cdot\id_A}\Hom(A,A)\rTo{\sim}\Tr_{\CC}(\Phi\circ\Psi)\to\Tr_{\CC}(F\circ\Phi\circ\Psi)\to \Tr_{\CC}(F).$$
Note that the last two arrows are obtained by applying $\Tr_{\CC}$ to the morphisms
$$A^\vee\boxtimes A\rTo{\id\ot\a} A^\vee\boxtimes F(A) \ \text{ and}$$ 
$$A^\vee\boxtimes F(A)\simeq F\circ (A^\vee\boxtimes A)\to F\circ \De_{\CC}.$$
Now the asserted equality follows from the fact that the last composition is equal to
\eqref{A-vee-F-A-eq}. 
\ed

Recall that we have an isomorphism 
\begin{equation}\label{A-FA-eq}
\Hom(A,F(A))\simeq \Hom(A^\vee, F'(A^\vee)),
\end{equation}
where $F'$ is given by the kernel $K'=K\circ\si$ (see Lemma \ref{F-op-lem}(i)).

\begin{lem}\label{tau-lem} 
Let $\a':A^\vee\to F'(A^\vee)$ be the morphism associated with a
morphism $\a:A\to F(A)$ under \eqref{A-FA-eq}. Then
$$\tau^A_F(\a)=\tau^{A^\vee}_{F'}(\a')$$
in $\Tr_{\CC}(F)\simeq\Tr_{\CC^{op}}(F')$.
\end{lem}

\Pf . The isomorphism \eqref{A-FA-eq} can be obtained from the isomorphisms
$$\Hom(A,F(A))\simeq\Tr_{\CC}((\Id_{\CC^{op}}\Box F)(A^\vee\boxtimes A))\simeq \Tr(F\circ X) \ \text{ and}$$
$$\Hom(A^\vee,F'(A^\vee))\simeq\Tr_{\CC^{op}}((\Id_{\CC}\Box F')(A\boxtimes A^\vee))\simeq
\Tr_{\CC}((F'\Box\Id_{\CC^{op}})(A^\vee\boxtimes A))\simeq \Tr(X\circ F),$$
where $X=A^\vee\boxtimes A$, together with the isomorphism $\Tr(X\circ F)\simeq\Tr(F\circ X)$.
Thus, it remains to check that the maps $F\circ X\to F$ and $X\circ F\to F$ induced by the 
canonical morphism $X\to\De_{\CC}$, give rise to the same maps after applying $\Tr$.
This boils down to the fact that the isomorphism 
$$\Tr(F)\simeq\Tr(\De\circ F)\simeq\Tr(F\circ \De)\simeq \Tr(F)$$ 
is equal to the identity.
\ed


As in the case of the generalized HRR-theorem of \cite[Lem.\ 1.3.2]{PV-HRR}, the fixed point formula
will be deduced from the following compatibility of the boundary-bulk maps with dg-functors.
 
\begin{lem}\label{main-lem}
Let $(\CC_1,F_1)$ and $(\CC_2,F_2)$ be a pair of dg-categories equipped with endofunctors. Suppose we have a dg-functor $(\Phi,\phi):(\per_{dg}(\CC_1),F_1)\to(\per_{dg}(\CC_2),F_2)$. 
For an object $A\in\per(\CC_1)$ 
let us consider the composed map 
$$(\Phi,\phi)_*^A:\Hom(A,F_1A)\rTo{\Phi}\Hom(\Phi(A),\Phi(F_1A))\to
\Hom(\Phi(A),F_2\Phi(A)),$$
where the second arrow is induced by $\phi$.
Then the following diagram is commutative
\begin{diagram}
\Hom(A,F_1A) &\rTo{(\Phi,\phi)_*^A} &\Hom(\Phi(A),F_2\Phi(A))\\
\dTo{\tau^A_{F_1}}&&\dTo{\tau^{\Phi(A)}_{F_2}}\\
\Tr(F_1)&\rTo{(\Phi,\phi)_*}&\Tr(F_2)
\end{diagram}
where the bottom horizontal arrow is the map \eqref{functoriality-map}.
\end{lem}

\Pf . The idea is to use the compatibility of the generalized
functoriality with compositions (see Lemma \ref{fun-comp-lem}).
Namely, as we have seen above, an element $\a\in\Hom(A,F_1A)$ gives rise to a dg-functor
$$(\iota_A,\a):(\per_{dg}(k),\Id)\to (\per_{dg}(\CC_1),F_1)$$
such that
$$\tau^A_{F_1}(\a)=(\iota_A,\a)_*(1)$$
(see Lemma \ref{bb-funct-lem}). Let $\a'=(\Phi,\phi)_*^A(\a)$.
It follows directly from the definitions that the dg-functor
$$(\iota_{\Phi(A)},\a')(\per_{dg}(k),\Id)\to (\per_{dg}(\CC_2),F_2)$$
is isomorphic to the composition
$(\Phi,\phi)\circ (\iota_A,\a)$.
Hence, by Lemma \ref{fun-comp-lem}, we obtain
$$(\iota_{\Phi(A)},\a')_*(1)=(\Phi,\phi)_*(\iota_A,\a)_*(1)=(\Phi,\phi)_*(\tau^A_{F_1}(\a)).$$
It remains to apply Lemma \ref{bb-funct-lem} to see that the left-hand side is equal to
$\tau^{\Phi(A)}_{F_2}(\a')$.
\ed

\subsection{Canonical pairings}\label{pairings-sec}

We can also use the functoriality to define a canonical pairing 
\begin{equation}\label{Tr-F-pairing}
\lan\cdot,\cdot\ran_{F,G}:\Tr(F)\ot\Tr(G)\to k,
\end{equation}
where $F:\per(\CC)\to\per(\CC)$ is given by a kernel $K$, and $G:\per(\CC)\to\per(\CC)$
is the right adjoint functor to $F$ associated with
the kernel $K^T$ (see Sec.\ \ref{trace-sec}).
Recall that $F':\per(\CC^{op})\to\per(\CC^{op})$ denotes the functor associated with the kernel
$K'=K\circ\si\in\per(\CC\ot\CC^{op}).$
We can apply the generalized functoriality \eqref{functoriality-map} to the functor
$$\Tr^{dg}_{\CC}:\per_{dg}(\CC^{op}\ot\CC)\to\per_{dg}(k),$$
the identity endofunctor of $\per_{dg}(k)$, and
the endofunctor $F'\Box G$ of $\per_{dg}(\CC^{op}\ot\CC)$.
Note that we have a morphism 
$$t(F,G):\Tr^{dg}_{\CC}\circ (F'\Box G)\to \Tr^{dg}_{\CC}$$
(this is \eqref{trace-nat-transformation} in the case $\DD=\CC$).
Thus, we can extend $\Tr^{dg}_{\CC}$ to a dg-functor
\begin{equation}\label{Tr-pairing-morphism-eq}
(\Tr^{dg}_{\CC},t(F,G)):(\per_{dg}(\CC^{op}\ot\CC),F'\Box G)\to (\per_{dg}(k),\Id).
\end{equation}
Hence, by the construction of Section \ref{funct-sec}, we obtain a map
$$(\Tr^{dg}_{\CC},t(F,G))_*:\Tr_{\CC^{op}\ot\CC}(F'\Box G)\to \Tr_k(\De_k)=k.$$
Using the isomorphisms
$$\Tr_{\CC^{op}\ot\CC}(F'\Box G)\simeq \Tr_{\CC^{op}}(F')\ot\Tr_{\CC}(G)$$
and $\Tr_{\CC^{op}}(F')\simeq \Tr_{\CC}(F)$ we obtain the pairing  \eqref{Tr-F-pairing}.

Applying the above construction to the adjoint pair $(G',F')$ (see Lemma 
\ref{G'-F'-lem})
we similarly get a pairing
$$\lan\cdot,\cdot\ran_{G',F'}:\Tr(G')\ot\Tr(F')\to k.$$

\begin{lem}\label{pairing-sym-lem} 
For $x\in\Tr(F)\simeq \Tr(F')$ and $y\in\Tr(G)\simeq\Tr(G')$ one has
$$\lan y,x\ran_{G',F'}=\pm\lan x,y\ran_{F,G},$$
where the sign is given by the Koszul rule.
\end{lem}

\Pf . This follows easily from the commutativity of the diagram \eqref{Tr-F-G-comm-diag}.
\ed

Recall that we have a canonical map
$$\ga_{F,G}:\De_{\CC}\to (F'\Box G)(\De_{\CC})$$
(see \eqref{ga-F-G0})

\begin{lem}\label{Casimir-lem} 
(i) The element 
$$\tau_{F,G}:=\tau^{\De_{\CC}}_{F'\Box G}(\ga_{F,G})\in \Tr(F')\ot \Tr(G)$$
satisfies 
\begin{equation}\label{pairing-tau-eq}
(\lan\cdot,\cdot\ran_{G',F'}\ot\id)(x\ot\tau_{F,G})=x,
\end{equation}
for any $x\in\Tr(G')\simeq\Tr(G)$.

\noindent
(ii) The pairing \eqref{Tr-F-pairing} is perfect, and $\tau_{F,G}\in \Tr(F')\ot\Tr(G)\simeq\Tr(F)\ot\Tr(G)$ 
is the corresponding Casimir element.
\end{lem}

\Pf . (i) By Lemma \ref{bb-funct-lem}, we can view $\tau^{\De_{\CC}}_{F'\Box G}(\ga_{F,G})$
as the generalized functoriality map associated with the dg-functor
$$(\de,\ga_{F,G}):(\per_{dg}(k),\Id)\to(\per_{dg}(\CC^{op}\ot\CC), F'\Box G)$$ 
where $\de:\per_{dg}(k)\to\per_{dg}(\CC^{op}\ot\CC)$ sends $k$ to $\De_\CC$.
It is easy to see that the composition
\begin{equation}\label{diagonal-comp-eq}
\per_{dg}(\CC)\rTo{\Id_{\CC}\Box\de} \per_{dg}(\CC\ot\CC^{op}\ot\CC)\rTo{\Tr_{\CC^{op}}\Box\Id_{\CC}}
\per_{dg}(\CC)
\end{equation}
is isomorphic to the identity. Indeed, for any kernel $K\in\per_{dg}(\CC^{op}\ot\CC)$
and $M\in\per_{dg}(\CC)$ we have
$$(\Tr_{\CC^{op}}\Box\Id_{\CC})(M\boxtimes K)\simeq M\otimes_\CC K.$$
Thus, for $K=\De_{\CC}$ we get
$$(\Tr_{\CC^{op}}\Box\Id_{\CC})(M\boxtimes\De_{\CC})\simeq M\otimes_\CC \De_{\CC}\simeq M.$$

Next, we observe that 
both arrows in \eqref{diagonal-comp-eq} extend to dg-functors between dg-categories with
endofunctors. Namely, we have dg-functors
$$(\Id_{\CC}\Box\de, \id\Box\ga_{F,G}): (\per_{dg}(\CC),G)\to (\per_{dg}(\CC\ot\CC^{op}\ot\CC),
G\Box F'\Box G)$$
and
$$(\Tr^{dg}_{\CC^{op}}\Box\Id_{\CC},t(G',F')\Box\id):
(\per_{dg}(\CC\ot\CC^{op}\ot\CC),G\Box F'\Box G)\to(\per_{dg}(\CC),G)$$
(where we use \eqref{Tr-pairing-morphism-eq} for the pair $(G',F')$).
We claim that the composition of these dg-functors is
\begin{equation}
(\Tr^{dg}_{\CC^{op}}\Box\Id_{\CC},t(G',F')\Box\id)\circ (\Id_{\CC}\Box\de, \id\Box\ga_{F,G})\simeq
(\Id_{\CC},\id).
\end{equation}
For this we have to check that the composition
\begin{align*}
&G(M)\simeq (\Tr_{\CC^{op}}\Box\Id_{\CC})(G(M)\boxtimes\De_{\CC})\rTo{\id\ot\ga_{F,G}}
(\Tr_{\CC^{op}}\Box\Id_{\CC})\circ(G\Box F'\Box G)(M\boxtimes\De_{\CC}) \\
&\rTo{t(G',F')\ot\id}(\Tr_{\CC^{op}}\Box G)(M\boxtimes\De_{\CC})\simeq 
G\circ (\Tr_{\CC^{op}}\Box \Id)(M\boxtimes\De_{\CC})\simeq G(M)
\end{align*}
is the identity. The corresponding composition of maps between the kernels can be rewritten as
$$K^T\simeq \Tr^{23}(K^T\boxtimes\De_\CC) \rTo{} \Tr^{23}(K^T\boxtimes (K\ot_\CC K^T))
\simeq \Tr^{23}(K^T\boxtimes K)\ot_\CC K^T\rTo{}\De_\CC\ot_\CC K^T\simeq K^T,$$
where $\Tr^{23}$ is the functor $\Tr_{\CC^{op}}$ applied in the 2nd and 3rd factors of the tensor product.
But the latter composition can be identified with the composition of two adjunction maps $G\to G\circ F\circ G\to 
G$, which is equal to the identity.

\noindent
(ii) The relation \eqref{pairing-tau-eq} shows that the left kernel of 
$\lan\cdot,\cdot\ran_{G',F'}$ is trivial. On the other hand,
applying (i) to the pair $(G',F')$ we obtain
$$(\lan\cdot,\cdot\ran_{F,G}\ot\id)(y\ot\tau_{G',F'})=y$$
for any $y\in\Tr(F)\simeq\Tr(F')$. Hence, the left kernel of $\lan\cdot,\cdot\ran_{F,G}$
is trivial. It remains to apply Lemma \ref{pairing-sym-lem}.
\ed

\begin{rems}\label{pairing-rems}
1. The above proof also shows that the tensor $\tau_{G',F'}$ is obtained from $\tau_{F,G}$
by the permutation of factors in the tensor product. This can be also seen directly from the definition.

\noindent
2. In the case $F=G=\Id_{\CC}$ both $\Tr(F)$ and $\Tr(G)$ are identified with the Hochschild homology $HH_*(\CC)$. Then the above canonical pairing gets identified with 
the pairing defined in \cite[Sec.\ 1.2]{PV-HRR}
(this pairing was introduced in \cite{Shk}, but we use the reversed order of argements)
using the canonical isomorphism 
\begin{equation}\label{HH-opposite-eq}
HH_*(\CC^{op})\simeq HH_*(\CC)
\end{equation}
(see \cite[Sec.\ 1.1]{PV-HRR}, \cite[Prop.\ 4.6]{Shk}).
\end{rems}


\subsection{Functoriality via inverse Serre functor}\label{Serre-sec}

The {\it inverse dualizing complex} of $\CC$ is
the object 
$$\De_{\CC^{op}}^{\vee}\in\per(\CC^{op}\ot\CC)$$
It is well known that this kernel induces the inverse of the Serre functor on $\per(\CC)$ (cf. \cite{Ginz-CY} or
\cite{Shkl-SD}). 

Recall that for any kernel $K\in\per(\CC^{op}\ot\CC)$ we have a natural isomorphism
\begin{equation}\label{Serre-Tr-eq}
\Hom(\De_{\CC^{op}}^{\vee},K)\simeq \Tr_{\CC}(K)
\end{equation}
(see \eqref{Serre-Tr0}).
We want to rewrite the definition of the functoriality map \eqref{functoriality-map} in terms of this isomorphism.

Thus, we assume that we have a dg-functor
$$(\Phi,\phi):(\per_{dg}(\CC_1),F_1)\to (\per_{dg}(\CC_2),F_2),$$
where $\phi:\Phi\circ F_1\to F_2\circ\Phi$ (see Section \ref{funct-sec}).
First, recall that there is a canonical morphism 
\begin{equation}\label{Serre-morphism}
\De_{\CC_2^{op}}^{\vee}\to (\Psi'\Box\Phi)(\De_{\CC_1^{op}}^{\vee})\simeq
\Phi\circ \De_{\CC_1^{op}}^{\vee}\circ\Psi
\end{equation}
(see \eqref{Serre-morphism0}).
Note also that the map $\phi$ gives rise to a map
\begin{equation}
\wt{\phi}: (\Psi'\Box\Phi)(F_1)=\Phi\circ F_1\circ\Psi\to F_2
\end{equation}
given as the composition 
$$\Phi\circ F_1\circ\Psi\rTo{\phi\circ\Psi} F_2\circ\Phi\circ\Psi\to F_2.$$

Now using isomorphism \eqref{Serre-Tr-eq}
we can give an alternative map $\Tr_{\CC_1}(F_1)\to \Tr_{\CC_2}(F_2)$
as the composition
\begin{equation}\label{Serre-fun-map}
\Hom(\De_{\CC_1^{op}}^{\vee},F_1)\to
\Hom((\Psi'\Box\Phi)(\De_{\CC_1^{op}}^{\vee}),(\Psi'\Box\Phi)(F_1))\to
\Hom(\De_{\CC_2^{op}}^{\vee},F_2),
\end{equation}
where the last arrow is induced by \eqref{Serre-morphism} and by $\wt{\phi}$.

\begin{prop}\label{Serre-prop} 
Under the isomorphism \eqref{Serre-Tr-eq} the map
\eqref{Serre-fun-map} gets identified with the map \eqref{functoriality-map}.
\end{prop}

\Pf . The map
$$\Hom(\De_{\CC_2^{op}}^{\vee},(\Psi'\Box\Phi)(F_1))\to\Hom(\De_{\CC_2^{op}}^{\vee},F_2)$$
induced by $\wt{\phi}$ gets identified under the isomorphism \eqref{Serre-Tr-eq} with
the composition
$$\Tr_{\CC_2}(\Phi\circ F_1\circ\Psi)\to\Tr_{\CC_2}(F_2\circ\Phi\circ\Psi)\to\Tr_{\CC_2}(F_2)$$
appearing as part of the definition of \eqref{functoriality-map}. 
Hence, it is enough to prove the commutativity of the
diagram
\begin{equation}\label{Serre-diagram}
\begin{diagram}
\Hom(\De_{\CC_1^{op}}^{\vee},F_1)&\rTo{}&
\Hom((\Psi'\Box\Phi)(\De_{\CC_1^{op}}^{\vee}),(\Psi'\Box\Phi)(F_1))&\rTo{}&
\Hom(\De_{\CC_2^{op}}^{\vee},(\Psi'\Box\Phi)(F_1))\\
\dTo{}&&&&\dTo{}\\
\Tr_{\CC_1}(F_1)&\rTo{}&\Tr_{\CC_1}(F_1\circ\Psi\circ\Phi)&\rTo{}&\Tr_{\CC_2}(\Phi\circ F_1\circ\Psi)
\end{diagram}
\end{equation}
Using the adjointness of the pair $(\Psi'\Box\Phi,\Phi'\Box\Psi)$ 
the first arrow in the top row can be identified with the map
$$\Hom(\De_{\CC_1^{op}}^{\vee},F_1)\to
\Hom(\De_{\CC_1^{op}}^{\vee},(\Phi'\Psi'\Box\Psi\Phi)(F_1)).$$
Next, by Lemma \ref{t-F-G-lem}(ii) the second arrow in the top row is equal to the following composition
(where we denoted $X=(\Psi'\Box\Phi)(F_1)$):
$$\Hom((\Psi'\Box\Phi)(\De_{\CC_1^{op}}^{\vee}),X)\to
\Hom(\De_{\CC_1^{op}}^{\vee},(\Phi'\Box\Psi)(X))\rTo{t(\Phi,\Psi)(X)}\Hom(\De_{\CC_1^{op}}^{\vee},X).$$
Therefore, by Lemma \ref{t-F-G-lem}(i), the composition of arrows in the top row of \eqref{Serre-diagram}
can be identified with the composition
$$\Tr(F_1)\to \Tr((\Phi'\Psi'\Box\Psi\Phi)(F_1))\rTo{\sim}\Tr(\Psi\Phi F_1\Psi\Phi)\rTo{\sim}
\Tr(\Phi\Psi\Phi F_1\Psi)\to\Tr(\Phi F_1\Psi).$$
Note that the composition of the first two arrows is equal to 
the composition of the following two maps induced by adjunction:
$$\Tr(F_1)\to \Tr(F_1\circ\Psi\circ\Phi)\to\Tr(\Psi\Phi F_1\Psi\Phi).$$
Hence, the commutativity of \eqref{Serre-diagram} follows from the commutativity of the outer
trapezoid in the diagram
\begin{diagram}
\Tr(\Psi\Phi F_1\Psi\Phi)&\rTo{}& \Tr(\Phi\Psi\Phi F_1\Psi)\\
\uTo{}&&\uTo{}&\rdTo{}\\
\Tr(F_1\Psi\Phi)&\rTo{}&\Tr(\Phi F_1\Psi)&\rTo{=}&\Tr(\Phi F_1\Psi)
\end{diagram}
in which the vertical arrows are induced by adjunction and the horizontal arrows in the left square are
isomorphisms of the form $\Tr(X\circ\Phi)\rTo{\sim}\Tr(\Phi\circ X)$. The
commutativity of the triangle in this diagram
follows from the fact that the composition $\Phi\to\Phi\Psi\Phi\to\Phi$ of
the maps induced by adjunction is equal to the identity.
\ed

In the situation when $F_1$ and $F_2$ are the identity functors we derive the following

\begin{cor} Let $\Phi:\per_{dg}(\CC_1)\to\per_{dg}(\CC_2)$ be a dg-functor.
Then the induced map 
$$\Phi_*:HH_*(\CC_1)\to HH_*(\CC_2)$$ 
defined as in Section \ref{funct-sec}, is equal to the composition
$$\Hom(\De_{\CC_1^{op}}^{\vee},\De_{\CC_1})\to
\Hom(\Phi\circ \De_{\CC_1^{op}}^{\vee}\circ\Psi, \Phi\circ\Psi)\to\Hom(\De_{\CC_2^{op}}^{\vee},\De_{\CC_2}),
$$
where the second arrow is induced by \eqref{Serre-morphism} 
and by the adjunction map $\Phi\circ\Psi\to \De_{\CC_2}$.
\end{cor}

\begin{rem}
One can deduce from the above corollary that our definition of functoriality maps on Hochschild homology 
is equivalent to the one given in \cite{CW} and \cite[Appendix]{BFK}.
\end{rem}

Recall that {\it Hochschild cohomology} of $\CC$ is defined by
$$HH^*(\CC)=\Hom_{\CC\ot\CC^{op}}(\De_{\CC},\De_{\CC}).$$
Note that for any autoequivalence 
$F:\per(\CC)\to\per(\CC)$ (given by a kernel) we 
we have a natural isomorphism
$$HH^*(\CC)\rTo{\sim} \Hom(F,F)$$
sending $\a$ to $F\circ\a$. Applying this to $F=\De_{\CC^{op}}^\vee$ we obtain an isomorphism
\begin{equation}\label{H-coh-isom}
HH^*(\CC)\simeq \Hom(\De_{\CC^{op}}^\vee,\De_{\CC^{op}}^\vee)\simeq\Tr_{\CC}(\De_{\CC^{op}}^\vee).
\end{equation}

Now if $\Phi:\per(\CC_1)\to\per(\CC_2)$ is an equivalence (given by a kernel) then
there is an induced isomorphism
\begin{equation}\label{Hoch-coh-fun}
\Phi_*:HH^*(\CC_1)\to HH^*(\CC_2).
\end{equation}
We are going to show that this map appears as one of our functoriality maps.
Since $\Phi$ is an equivalence, we have a natural isomorphism
$$\wt{\phi}:\Phi\circ\De_{\CC_1^{op}}^\vee\circ\Phi^{-1}\rTo{\sim}\De_{\CC_2^{op}}^\vee,$$
and hence the induced isomorphism
\begin{equation}\label{H-coh-Serre-eq}
\phi:\Phi\circ\De_{\CC_1^{op}}^\vee\rTo{\sim} \De_{\CC_2^{op}}^\vee\circ\Phi.
\end{equation}

\begin{cor}\label{H-coh-cor} 
The map \eqref{Hoch-coh-fun} gets identified under the isomorphisms \eqref{H-coh-isom}
with
$$(\Phi,\phi)_*:\Tr_{\CC_1}(\De_{\CC_1^{op}}^\vee)\to\Tr_{\CC_2}(\De_{\CC_2^{op}}^\vee).$$
\end{cor}

\Pf . By Proposition \ref{Serre-prop}, it is enough to prove the commutativity of
the diagram
\begin{diagram}
\Hom(\De_{\CC_1},\De_{\CC_1})&\rTo{\a\mapsto \Phi\circ\a\circ\Phi^{-1}}&\Hom(\De_{\CC_2},\De_{\CC_2})\\
\dTo{\sim}&&\dTo{\sim}\\
\Hom(\De_{\CC_1^{op}}^\vee,\De_{\CC_1^{op}}^\vee)&\to\Hom(\Phi\circ\De_{\CC_1^{op}}^\vee\circ\Phi^{-1},
\Phi\circ\De_{\CC_1^{op}}^\vee\circ\Phi^{-1})\to& \Hom(\De_{\CC_2^{op}}^\vee,\De_{\CC_2^{op}}^\vee)
\end{diagram}
where the vertical arrows are isomorphisms \eqref{H-coh-isom}
and the second arrow in the bottom row is induced by the isomorphism $\wt{\phi}$.
Upon unraveling the definitions this becomes a tautology.
\ed

In the remainder of this section we explain the relation between the canonical pairing on the Hochschild homology
and the Mukai pairing considered in \cite{Cal-I,CW,Ram2}. This will not be used anywhere else in
the paper. Let us recall the definition of the Mukai pairing $\lan\cdot,\cdot\ran_M$ on $HH_*(X)$,
where $X$ is a smooth projective variety, following \cite{Cal-I}\footnote{The definition of the Mukai pairing in \cite[Sec.\ 5]{CW} seems to contain a misprint}. Let $S_X$ (resp., $S^{-1}_X$) 
be the kernel in $D^b(X\times X)$ corresponding to the
Serre functor $\SS_X$ (resp., the inverse of the Serre functor $\SS^{-1}_X$) on $D^b(X)$. We have natural isomorphisms
$$HH_*(X)\simeq \Hom^*_{X\times X}(S^{-1}_X,\De_X),$$
$$\SS_{X\times X}(S^{-1}_X)\simeq S_X,$$
where $\De_X$ is the structure sheaf of the diagonal in $X\times X$ (see \cite[Sec.\ 2.2, Prop.\ 1]{CW}).
From the second isomorphism we get using Serre duality on $X\times X$ a canonical functional
$$\tr_X:\Hom_{X\times X}(S^{-1}_X,S_X)\to k.$$
Now the Mukai pairing is given by
$$\lan x,y\ran_M:= \tr_X(\tau(x)\circ y),$$
where $x,y\in\Hom_{X\times X}(S^{-1}_X,\De_X)$, and $\tau$ ($=\tau_L$ in the notation of \cite{CW})
is the isomorphism 
$$\tau:\Hom_{X\times X}(S^{-1}_X,\De_X)\to \Hom_{X\times X}(\De_X,S_X):x\mapsto x\circ S_X.$$
On the other hand, if we take $\CC$ to be a dg-enhancement of $D^b(X)$, we get an identification
$HH_*(\CC)\simeq HH_*(X)$, so we can compare the Mukai pairing with 
the canonical pairing $\lan\cdot,\cdot\ran$ on $HH_*(\CC)$
(see Remark \ref{pairing-rems}.2).

\begin{prop} With the above identifications we have the equality of pairings
$\lan\cdot,\cdot\ran_M=\lan\cdot,\cdot\ran$ on $HH_*(X)$.
\end{prop}

\Pf . By definition, $\lan\cdot,\cdot\ran$ is obtained by looking at the map 
\begin{equation}\label{pre-pairing-map}
\Hom_{X^4}(S^{-1}_{X^2},\De_{X^2})=HH_*(X^2)\to k
\end{equation}
on Hochschild homology induced by the functor
$$\Tr_X=R\Hom_{X^2}(S_X^{-1},?):D^b(X\times X)\to D^b(k).$$
The right adjoint to this functor is given by $S_X\in D^b(X\times X)$, 
and the description of the functoriality via Serre functors \eqref{Serre-fun-map} 
implies that \eqref{pre-pairing-map} is given by
$$\a\mapsto \tr_X((\Tr_X\circ\a)_{S_X}(u)),$$ 
where $\a\in\Hom(S^{-1}_{X^2},\De_{X^2})$, the canonical element
$$u\in\Tr_X(\SS^{-1}_{X^2}(S_X))\simeq\Hom_{X^2}(S_X^{-1},S_X^{-1})$$
corresponds to the identity map of $S^{-1}_X$, and
$$(\Tr_X\circ\a)_{S_X}:(\Tr_X\circ \SS^{-1}_{X^2}) (S_X)\to \Tr_X(S_X)$$
is the map induced by the morphism of functors $\Tr_X\circ\a$.
Note that the action of $\SS^{-1}_{X^2}$ on kernels in $D^b(X\times X)$ is given by
$$\SS^{-1}_{X^2}(F)=S^{-1}_X\circ F\circ S^{-1}_X.$$
Thus, $(\Tr_X\circ\a)_{S_X}(u)$ is the composition
$$S_X^{-1}\simeq S^{-1}_X\circ S_X\circ S^{-1}_X\rTo{\a_{S_X}} S_X.$$
Finally, the K\"unneth isomorphism 
$$\Hom_{X^2}(S^{-1}_X,\De_X)^{\ot 2}\to \Hom_{X^4}(S^{-1}_{X^2},\De_{X^2})$$
sends $x\ot y$ to the morphism $\a:S^{-1}_{X^2}\to \De_{X^2}$ such that
$\a_F$ is given by the composition
$$S^{-1}_X\circ F\circ S^{-1}_X\rTo{S^{-1}_X\circ F\circ y} S^{-1}_X\circ F\rTo{x\circ F} F,$$ 
which in the case $F=S_X$ gets identified with the composition
$$S^{-1}_X\rTo{y} \De_X\rTo{x\circ S_X} S_X.$$
Applying the map $\tr_X$ we get exactly the Mukai pairing $\lan x,y\ran_M$.
\ed

\begin{rems} 1. Recall that the categorical version of the Hirzebruch-Riemann-Roch formula states
that for a pair of objects $A,B\in\per(\CC)$ one has
$$\lan \ch(A^\vee),\ch(B)\ran=\chi(\Hom(A,B))$$
(see \cite[(1.2)]{Shk}, \cite[(1.18)]{PV-HRR}).
In \cite[Thm.\ 14]{CW} 
the left-hand side is replaced with $\lan \ch(A),\ch(B) \ran_M$. The fact that the two formulas
are the same follows from the equality
\begin{equation}\label{Chern-dual-eq}
\ch(A)=\ch(A^\vee),
\end{equation}
given by Lemma \ref{tau-lem} in the case $F=\Id_{\CC}$, $\a=\id_A$.

\noindent
2. In the case when $\CC$ is the dg-enhancement of $D^b(X)$, where $X$ is smooth and projective,
there is an equivalence $\CC^{op}\simeq\CC$ given by the duality functor 
$\FF\mapsto R\underline{Hom}(\FF,\OO_X)$. This equivalence induces an isomorphism
$HH_*(\CC^{op})\simeq HH_*(\CC)$, which is different from the abstract isomorphism. This is the reason
for the extra automorphism $\vee$ of $HH_*(X)$ appearing in \cite[Thm.\ 1]{Ram2}.
The same automorphism appears if one converts \eqref{Chern-dual-eq} into a formula
for the Chern character of the dual vector bundle on $X$.
\end{rems}

\section{Lefschetz type formulas}

\subsection{Categorical version of holomorphic Lefschetz formula}

Now let $A$ and $B$ be a pair of objects of $\per(\CC)$ equipped with morphisms $\alpha:A\to F(A)$ and $\beta:F(B)\to B$. Let $G$ be the right adjoint to $F$, and let $F'$ be the functor associated with the kernel $K':=K\circ\si$ (see Section \ref{adj-sec}).

\begin{thm} In this situation the formula \eqref{HLF-eq} holds, where 
$\tau^A_F$ (resp., $\tau^B_G$) is the generalized boundary
bulk map of Section \ref{bbm-sec}, and $\lan\cdot,\cdot\ran_{F,G}$ is the canonical pairing of
Section \ref{pairings-sec}.
\end{thm}

\Pf .
We would like to apply Lemma \ref{main-lem} to the dg-functor
\eqref{Tr-pairing-morphism-eq}, extending $\Tr^{dg}_{\CC}$,
from $(\per_{dg}(\CC^{op}\ot\CC), F'\Box G)$ to $(\per_{dg}(k),\Id)$.

By Lemma \ref{F-op-lem}, the map $\alpha:A\to F(A)$ in $\per(\CC)$ induces a morphism 
$\alpha':A^\vee\to F'(A^\vee)$ in $\per(\CC^{op})$.
Also, $\b:F(B)\to B$ induces a map $\wt{\b}:B\to G(B)$.
Hence, we obtain a morphism
$$\alpha'\ot\wt{\beta}:A^\vee\boxtimes B\to (F'\Box G)(A^\vee\boxtimes B).$$
Note that
$$\tau^{A^\vee\boxtimes B}_{F'\Box G}=\tau^{A^\vee}_{F'}\ot \tau^B_G.$$
Also, for a graded vector space $V$ the map
$$\tau^V_{\Id}:\Hom(V,V)\to k$$
is the usual supertrace.
Hence, Lemma \ref{main-lem} gives an equality 
\begin{equation}\label{trace-prelim-eq}
\sTr(f)=\lan \tau^{A^\vee}_{F'}(\a'),\tau^B_G(\wt{\b})\ran_{F,G},
\end{equation}
where $f=(\Tr^{dg}_{\CC},t(F,G))^{A^\vee\boxtimes B}_*(\a'\ot\wt{\b})$ is the composition
$$\Hom(A,B)\simeq
\Tr_{\CC}(A^\vee\boxtimes B)\rTo{\Tr_{\CC}(\a'\ot\wt{\b})} 
\Tr_{\CC}(F'(A^\vee)\boxtimes G(B))\rTo{t(F,G)}\Tr_{\CC}(A^\vee\boxtimes B)\simeq\Hom(A,B).$$
Since $\tau^{A^\vee}_{F'}(\a')=\tau^A_F(\a)$ by Lemma \ref{tau-lem}, the right-hand side
of \eqref{trace-prelim-eq} coincides with the right-hand side of \eqref{HLF-eq}.
It remains to check that $f$ is equal to the endomorphism 
$(F,\a,\b)_*$ of $\Hom(A,B)$ (see \eqref{F-a-b-eq}). We can rewrite the definition of $(F,\a,\b)_*$ as
the following composition:
$$\Hom(A,B)\simeq
\Tr_{\CC}(A^\vee\boxtimes B)\rTo{s}
\Tr_{\CC}(F(A)^\vee\boxtimes F(B))\rTo{\Tr_{\CC}(\a^\vee\ot\b)}
\Tr_{\CC}(A^\vee\boxtimes B)\simeq\Hom(A,B),$$
where $s$ is induced by the isomorphism $F(A)^\vee\simeq G'(A^\vee)$ (see \eqref{G'-F-dual-eq}) and
by a natural morphism of functors
$$s(G,F):\Tr_{\CC}\to\Tr_{\CC}\circ(G'\Box F).$$
This morphism of functors is constructed similarly to $t(F,G)$ using the map of
kernels \eqref{adj-from-id}.
Now the required equality $f=(F,\a,\b)_*$ follows from the commutativity of the diagram
\begin{diagram}
\Tr_{\CC}(A^\vee\boxtimes B)&\rTo{\Tr(\a'\ot\wt{\b})}&\Tr_{\CC}(F'(A^\vee)\boxtimes G(B))\\
\dTo{s(G,F)}&&\dTo{s(G,F)}\\
\Tr_{\CC}(G'(A^\vee)\boxtimes F(B))&\rTo{\Tr(G'(\a')\ot F(\wt{\b}))}&
\Tr_{\CC}(G'F'(A^\vee)\boxtimes FG(B))\\
\dTo{}&&\dTo{}\\
\Tr_{\CC}(F(A)^\vee\boxtimes F(B))&\rTo{\Tr_{\CC}(\a^\vee\ot\b)}&\Tr_{\CC}(A^\vee\boxtimes B)
\end{diagram}
where the vertical arrows in the lower square are induced by the isomorphism
$G'(A^\vee)\simeq F(A)^\vee$ and by the adjunction morphisms
$G'F'\to\Id$ and $FG\to\Id$.
Note that the commutativity of the lower square follows from Lemma \ref{F-op-lem}(ii).
\ed



\begin{ex}\label{classical-ex} 
The classical setup is recovered if we take $\CC$ to be a dg-enhancement of
the derived category $D^b(M)$ of coherent sheaves on a smooth projective variety $M$ and
$F=f^*$, where $f:M\to M$ is an endomorphism. The category of kernels in this case can
be identified with $D^b(M\times M)$, so that the functor $\Tr:D^b(M\times M)\to \per(k)$
sends $K$ to $H^*(M,\De^*(K))$. In particular, the space
$$T_f:=\Tr(f^*)=\Tr(f_*)=H^*(M,\De^*(\OO_{\Ga_f}))$$
corresponds to considering the derived intersection of the graph $\Ga_f$ with the
diagonal. In the case when the intersection is transversal we have 
$$T_f=\bigoplus_{f(x)=x} k.$$
If $V$ is a vector bundle (or a complex of vector bundles) on $M$ then applying $(f\times\id_M)^*$
to the natural map $(V^\vee\boxtimes V)\to \OO_\De$ and then restricting to the diagonal
we get a map
$$\tau_f^V:\Hom(V,f_*V)=\Hom(f^*V,V)\to T_f$$
which is exactly the map used in \eqref{HLF-eq}.
For example, for $V=\OO_X$ we get $\tau_f^{\OO_X}(\id)=\unit\in T_f$, the element obtained from
the projection $\OO_X\to \OO_{\Ga_f}$. In the case of transversal fixed points we have
$$\tau_f^V(\a)=\sum_{f(x)=x} \tr(\a, V|_x)\cdot \de_x,$$
where $(\de_x)$ is the natural basis of $T_f$ indexed by fixed points. Finally, we claim that in
the case of transversal intersection the canonical form $\lan\cdot,\cdot\ran$ on $T_f$ is
given by 
$$\lan \de_x,\de_x\ran=\frac{1}{\det(\id-d_xf)} \ \text{ and }\ \lan\de_x,\de_y\ran=0 \ \text{ for }\ x\neq y,$$
where $d_xf: T_xM\to T_xM$ is the tangent map to $f$ at $x$.
Indeed, the fact that the basis $(\de_x)$ is orthogonal can be deduced from 
Lemma \ref{Casimir-lem}. To get the formula for $\lan \de_x,\de_x\ran$, where $f(x)=x$, we observe that
$\tau_f^{\OO_x}=c_x\cdot\de_x$ for some $c_x\in k$, and 
apply the formula \eqref{HLF-eq} to the pairs $(\OO_M, \OO_x)$ and to $(\OO_x,\OO_x)$. 
In the former case we get
$$1=\lan \unit, c_x\de_x\ran=c_x\cdot\lan\de_x,\de_x\ran,$$
and in the latter case we get
$$\sTr(F_*,\Ext^*(\OO_x,\OO_x))=c_x^2\cdot\lan \de_x,\de_x\ran.$$
Hence, 
$$c_x=\lan \de_x,\de_x\ran^{-1}=\sTr(f_*,\Ext^*(\OO_x,\OO_x)).$$
But $\Ext^*(\OO_x,\OO_x)$ is the exterior algebra of $\Ext^1(\OO_x,\OO_x)\simeq T_xM$
and the action of $F$ on $\Ext^1(\OO_x,\OO_x)$ is precisely the tangent map $d_xf$.
Since $F_*$ is compatible with the algebra structure, we obtain
$$\sTr(F_*,\Ext^*(\OO_x,\OO_x))=\det(\id-d_xf).$$
Thus, our formula in this case becomes the classical holomorphic Lefschetz formula
$$\sTr(f,H^*(M,V))=\sum_{f(x)=x} \tr(f,V|_p)\cdot\det(\id-d_xf)^{-1}.$$
\end{ex}

Using the Lefschetz formula \eqref{HLF-eq} we can give some criteria for non-vanishing of
the space $\Tr_{\CC}(F)$ (which is roughly analogous to the existence of a fixed point of
an endomorphism).

Recall that an object $A$ of a $k$-linear
triangulated category is called {\it exceptional} if $\Hom^i(A,A)=0$ for $i\neq 0$ and
$\Hom^0(A,A)=k$. An object $A$ is called {\it $n$-spherical} if $\Hom^i(A,A)=0$ for $i\neq 0,n$
and $\Hom^0(A,A)=\Hom^n(A,A)=k$.

\begin{cor} One has $\Tr_{\CC}(F)\neq 0$ if one of the following conditions holds:

\noindent
(i) there exists an exceptional object $A\in\per(\CC)$ such that $F(A)\simeq A$;

\noindent
(ii) there exists an $n$-spherical object $A\in\per(\CC)$, where $n$ is even, such that $F(A)\simeq A$,
and in addition $F^N\simeq\Id$, where $N$ is odd, and the characteristic of $k$ is not $2$;

\noindent
(iii) there exists an object $A\in\per(\CC)$ with $\Hom(A,A)\simeq\bigwedge^*_k(\Hom^1(A,A))$,
such that $F(A)\simeq A$ and the endomorphism $F_*$ of $\Hom^1(A,A)$ satisfies
$\det(\id-F_*)\neq 0$.
\end{cor}

\Pf . We apply \eqref{HLF-eq} to $A=B$, with $\a:A\to F(A)$ some isomorphism and $\b=\a^{-1}$.
We have to show that in all three cases we have $\sTr(F_*,\Hom(A,A))\neq 0$.

In case (i) we have $\Hom(A,A)=k$ and the induced map $F_*$ on it 
is identity, so $\sTr(F_*,\Hom(A,A))=1$.

In case (ii) $F_*$ acts as identity on $\Hom^0(A,A)=k\cdot\id_A$ 
and as some $N$th root of unity $\zeta$ on $\Hom^n(A,A)\simeq k$. Since $n$ is even,
we obtain 
$$\sTr(F_*,\Hom(A,A))=1+\zeta$$
which is nonzero since $N$ is odd.

Finally, in case (iii) we use the fact that $F_*$ is an endomorphism of the algebra
$\Hom(A,A)\simeq\bigwedge^*_k(\Hom^1(A,A))$, so we get
$$\sTr(F_*,\Hom(A,A))=\det(\id-F_*)\neq 0.$$
\ed

Note that the property (iii) in the above Corollary is modeled on the properties of
the pair $(f^*,\OO_x)$ for a transversal fixed point $x$ of an endomorphism $f$.

\subsection{Lefschetz reciprocity}\label{top-sec}

As before, we assume that $\CC$ is a dg-category such that $\per_{dg}(\CC)$ is saturated,
and $F$ is a dg-endofunctor of $\per_{dg}(\CC)$.
Now let $\Psi$ be another dg-endofunctor of $\per_{dg}(\CC)$ together with a morphism
$$f: F\circ\Psi\to \Psi\circ F$$
Then we have the induced endomorphism
$(F,f)_*$ of $\Tr(\Psi)$. Let $G$ be right adjoint to $F$. Then we have the induced map
$$\psi: \Psi\circ G\to G\circ F\circ\Psi\circ G\rTo{G\circ f\circ G} G\circ\Psi\circ F\circ G\to  G\circ\Psi,$$
and hence the induced endomorphism $(\Psi,\psi)_*$ of $\Tr(G)$.
We have the following {\it Lefschetz reciprocity} formula.

\begin{thm}\label{LR-thm} In the above situation
\begin{equation}\label{LR-eq}
\sTr((F,f)_*,\Tr(\Psi))=\sTr((\Psi,\psi)_*,\Tr(G)).
\end{equation}
\end{thm}

\Pf .
The idea is to apply formula \eqref{HLF-eq} to the category $\CC^{op}\ot\CC$,
equipped with the endofunctor $G'\Box F$,
taking as $B$ the kernel representing $\Psi$ and 
$$A=\De_{\CC^{op}}^{\vee},$$
the inverse dualizing complex (see section \ref{Serre-sec}).

Note that we have
$$\Tr_{\CC^{op}\ot\CC}(G'\Box F)\simeq\Tr_{\CC^{op}}(G')\ot \Tr_{\CC}(F)\simeq
\Tr(G)\ot\Tr(F).$$
Also, the right adjoint to $G'\Box F$ is $F'\Box G$ (by Lemma \ref{G'-F'-lem}).

We will apply \eqref{HLF-eq} to certain natural morphisms 
$$\a:\De_{\CC^{op}}^\vee\to (G'\Box F)(\De_{\CC^{op}}^\vee) \text{ and}$$
$$\b:(G'\Box F)(B)\to B$$
that we are going to define presently.
Namely, using the isomorphism of Lemma \ref{F-op-lem}(i), $\a$ corresponds to a morphism
$$\a':\De_{\CC^{op}}\to (G'\Box F)'(\De_{\CC^{op}})=(G\Box F')(\De_{\CC^{op}}),$$
and we take 
$$\a':=\ga_{G',F'}$$
(see \eqref{ga-F-G0}).
On the other hand, $\b$ corresponds to the morphism of functors
$$F\circ \Psi\circ G\rTo{f\circ G} \Psi\circ F\circ G\to \Psi.$$
If follows that the map 
$$\wt{\b}:B\to (F'\Box G)(B),$$
obtained from $\b$ by adjunction of the functors $(G'\Box F, F'\Box G)$, corresponds to the morphism
of functors
$$\Psi\to \Psi\circ G\circ F\rTo{\psi\circ F} G\circ\Psi\circ F.$$


It follows from Proposition \eqref{Serre-prop} that the endomorphism $(G'\Box F,\a,\b)_*$ of 
$\Hom(A,B)=\Hom(\De_{\CC^{op}}^\vee,B)=\Tr(\Psi)$
can be identified with the endomorphism $(F,f)_*$ of $\Tr(\Psi)$.

Thus, the formula \eqref{HLF-eq} combined with Lemma \ref{tau-lem} gives
$$\sTr((F,f)_*,\Tr(\Psi))=
\lan \tau^{\De_{\CC^{op}}}_{G\Box F'}(\ga_{G',F'}), \tau^{B}_{F'\Box G}(\wt{\b})\ran=
\lan \tau_{G',F'}, \tau^{B}_{F'\Box G}(\wt{\b})\ran.$$

We claim that
\begin{equation}\label{2-tau-eq}
\tau^{B}_{F'\Box G}(\wt{\b})=(\id\ot (\Psi,\psi)_*)(\tau_{F,G}),
\end{equation}
where
$$\psi:\Psi\circ G\to G\circ\Psi$$
is induced by $f$.
Indeed, this follows from Lemma \ref{main-lem} since $B=(\Id\Box\Psi)(\De_{\CC})$. 

Applying Lemma \ref{str-lem} below we get
$$\lan \tau_{G',F'},(\id\ot (\Psi,\psi)_*)(\tau_{F,G})\ran=\sTr((\Psi,\psi)_*,\Tr(G))$$
It remains to observe that starting with the data
$(F,\Psi,f)$ we can recover the data $(\Phi,F,\phi)$. Hence,
we can consider $(F,\Psi,f)$ as the primary data. Now
\eqref{LR-eq} is obtained from the above equality by a change of notation.
\ed

\begin{lem}\label{str-lem} Let $V$ and $W$ be finite-dimensional
supervector spaces equipped with a perfect (even) pairing
$$\lan\cdot,\cdot\ran_{V,W}: V\ot W\to k.$$
Let $\tau\in V\ot W$ be the corresponding Casimir element, 
and let $\tau'=\si(\tau)\in W\ot V$. 
Consider the induced pairing
$$\lan\cdot,\cdot\ran_{V\ot W, W\ot V}: (V\ot W)\ot (W\ot V)\to k$$
given by 
$$\lan v\ot w, w'\ot v'\ran_{V\ot W, W\ot V}=\pm \lan v,w'\ran_{V,W}\cdot \lan v',w\ran_{V,W}$$
with the sign given by the Koszul rule.
Then for any even automorphism $A$ of $V$ one has
$$\lan \tau, (\id\ot A)(\tau')\ran_{V\ot W, W\ot V}=\sTr(A).$$
In particular,
$$\lan \tau, \tau'\ran_{V\ot W, W\ot V}=\sdim(V).$$
\end{lem}

\Pf . Let $(v_i)$ be a homogeneous basis in $V$ and let $(w_i)$ be the dual basis of $W$, so that
$\lan v_i,w_j\ran=\de_{ij}$.
Then 
$$\tau=\sum_i (-1)^{\deg(v_i)} v_i\ot w_i \ \text{ and } \ \tau'=\sum_i w_i\ot v_i.$$
It is enough to prove the formula for $A$ such that $A(v_i)=v_j$ and $A(v_l)=0$ for $l\neq i$.
Then we have 
$$\lan \tau, (\id\ot A)(\tau')\ran_{V\ot W, W\ot V}=
\lan (-1)^{\deg(v_i)}v_i\ot w_i, w_i\ot v_j\ran_{V\ot W, W\ot V}=\de_{ij}(-1)^{\deg(v_i)}=\sTr(A).$$
\ed

\begin{rem}
In the particular case $(\Psi,\psi)=(\Id_{\CC},\id)$ the formula \eqref{LR-eq} 
gives \eqref{LF-eq}.
Note however, that the argument of \cite{Lunts} for the proof of \eqref{LF-eq} is more direct.
\end{rem}

\begin{cor}\label{FNabs-cor} Keep the assumptions of Theorem \ref{LR-thm}. Set
$\wt{\Psi}=\Psi\circ F$ and let
$$\wt{f}:F\circ\wt{\Psi}=F\circ \Psi\circ F \to \Psi\circ F\circ F=\wt{\Psi}\circ F$$
be the map $f\circ F$. Then
\begin{equation}\label{FNabs-eq}
\sTr((F,f)_*,\Tr(\Psi))=\sTr((F,\wt{f}),\Tr(\wt{\Psi})).
\end{equation}
A similar result holds for $F\circ\Psi$ instead of $\Psi\circ F$.
\end{cor}

\Pf . By Theorem \ref{LR-thm}, the equality \eqref{FNabs-eq} is equivalent to
$$\sTr((\Psi,\psi)_*,\Tr(G))=\sTr((\wt{\Psi},\wt{\psi})_*,\Tr(G)),$$
where
$$\wt{\psi}:\wt{\Psi}\circ G\to G\circ\wt{\Psi}$$
is defined similarly to $\psi$, starting from $\wt{f}$.
We claim that in fact
\begin{equation}\label{Psi-*-G-eq}
(\Psi,\psi)_*=(\wt{\Psi},\wt{\psi})_*.
\end{equation}
Let us first calculate the morphism $\wt{\psi}$. We have a commutative diagram
\begin{diagram}
\Psi\circ F\circ G &\rTo{}& G\circ F\circ\Psi\circ F\circ G &\rTo{G\circ f\circ F\circ G}& G\circ \Psi\circ F\circ F\circ G\\
\dTo{}&&\dTo{}&&\dTo{}\\
\Psi&\rTo{}& G\circ F\circ \Psi&\rTo{G\circ f}& G\circ\Psi\circ F
\end{diagram}
Here and in the diagrams below unnamed arrows are induced by adjunction morphisms.
By definition, $\wt{\psi}$ is given by the composition of the arrows in the top row followed by the right vertical arrow. Hence, it can also be computed as the composition of the arrows in the bottom row with the left
vertical arrow. Now the commutative diagram
\begin{diagram}
\Psi &\rTo{}& G\circ F\circ \Psi&\rTo{G\circ f}& G\circ\Psi\circ F\\
\dTo{}&&\dTo{}&&\dTo{}&\rdTo{\id}\\
\Psi\circ G\circ F &\rTo{}&G\circ F\circ\Psi\circ G\circ F&\rTo{G\circ f\circ G\circ F}&
G\circ\Psi\circ F\circ G\circ F&\rTo{}&G\circ\Psi\circ F
\end{diagram}
together with the definition of $\psi$, shows that $\wt{\psi}$ is equal to the composition
$$\Psi\circ F\circ G\rTo{\Psi\circ c} \Psi\circ G\circ F\rTo{\psi\circ F} G\circ\Psi\circ F,$$
where
$$c: F\circ G\to \Id\to G\circ F$$
is the composition of two adjunction morphisms.
Hence, by Lemma \ref{fun-comp-lem}, \eqref{Psi-*-G-eq} would follow once we prove that the map
$$(F,c)_*:\Tr(G)\to\Tr(G)$$
is the identity map.
Recall (see Sec.\ \ref{funct-sec}) that the map $(F,c)_*$ is given by the composition
$$\Tr(G)\to\Tr(G\circ F\circ G)\rTo{\si}\Tr(F\circ G\circ G)\rTo{\Tr(c\circ G)}\Tr(G\circ F\circ G) \to \Tr(G),$$
Here and below arrows marked by $\si$ are induced by the isomorphism \eqref{Tr-switch}.
By definition of the morphism $c$ this is equal to the composition
\begin{equation}\label{Tr-G-comp}
\Tr(G)\to\Tr(G\circ F\circ G)\rTo{\si}\Tr(F\circ G\circ G)\to\Tr(G)\to\Tr(G\circ F\circ G) \to \Tr(G).
\end{equation}
Note that the composition of the last two arrows is the identity map, which together with the commutative triangle
\begin{diagram}
\Tr(G\circ F\circ  G)\\
\dTo{\si}&\rdTo{}\\
\Tr(F\circ  G\circ G)&\rTo{}&\Tr(G)
\end{diagram}
shows that the composition \eqref{Tr-G-comp} is equal to the identity, i.e., $(F,c)_*=\id$. 
\ed

Consider again the situation when $\CC$ is a dg-enhancement of $D^b(M)$, where $M$ is a smooth projective variety over $k$, $F=f^*$ and $\Psi(\FF)=V\ot_{\OO_M}\FF$, where $V$ is a bounded complex of vector bundles
on $M$. 
Thus, the kernel on $M\times M$ corresponding to $\Psi$ is $\De_*V$,
where $\De:V\to V\times V$ is the diagonal embedding. A map
$F\circ\Psi\to\Psi\circ F$ is the same as a map $\a:f^*V\to V$, and there is a natural isomorphism
$$\Tr(\Psi)\simeq H^*(M,\De^*\De_*V)\simeq H^*(M,\Om_M^*\ot V).$$
The map $(F,\a)_*:\Tr(\Psi)\to\Tr(\Psi)$ induced by $\a$ is the composition
$$\a_*:H^*(M,\Om_M^*\ot V)\rTo{f^*} H^*(M,\Om_M^*\ot f^*V)\rTo{\a} H^*(M,\Om_M\ot V).$$
On the other hand, the functor $F=f^*$ is given by the kernel
$(f,\id_M)_*\OO_M$ on $M\times M$, so
$\Psi\circ F$ is given by the kernel $(f,\id_M)_*V$.
Hence,
$$\Tr(\Psi\circ F)\simeq H^*(M,\De^*(f,\id_M)_*V)\simeq H^*(M,i_*i^*V)\simeq H^*(M^f,i^*V),
$$
where $M^f$ is the {\it derived intersection} of $\Ga_f$ with the diagonal in $M\times M$,
fitting into the derived cartesian diagram
\begin{diagram}
M^f &\rTo{i}& M\\
\dTo{i} &&\dTo{(\id_M,f)}\\
M&\rTo{\De}& M\times M
\end{diagram}
Corollary \ref{FNabs-cor} leads in this situation to the following result conjectured by Frenkel-Ng\^{o}
(see \cite[Conj.\ 6.2]{FN}).

\begin{thm}\label{der-int-thm}
Let $V$ be a bounded complex of vector bundles on a smooth projective variety $M$ over $k$,
$f:M\to M$ and endomorphism. Then for any moprhism $\a:f^*V\to V$ in $D^b(M)$ one has
\begin{equation}\label{der-int-formula}
\sTr(\a_*,H^*(M,\Om_M^*\ot V))=\sTr((i^*\a)_*:H^*(M^f,i^*V)),
\end{equation}
where the map $(i^*\a)_*$ is given by
$$H^*(M^f,i^*V)\simeq H^*(M^f,i^*f^*V)\rTo{i^*\a} H^*(M^f,i^*V)$$
\end{thm}


\begin{rems}
1. In the case of transversal intersection the formula \eqref{der-int-formula} becomes
$$\sTr(f,H^*(M,\Om_M^*\ot V))=\sum_{f(x)=x}\sTr(f,V|_p)$$
which follows also from the usual holomorphic Lefschetz formula (see Example \ref{classical-ex}).

\noindent
2. In the formula of \cite[Conj.\ 6.2]{FN} 
the right-hand side apparently has a typo: it should have $i^*V$ rather than
$i^*\De_*(V)$ in order to be correct even in the transversal case.

\noindent
3. Theorem \ref{der-int-thm} 
and its proof should generalize to an appropriate class of derived algebraic stacks, as envisioned in \cite{FN}.
\end{rems}

\subsection{Lefschetz formula in Hochschild cohomology}\label{H-coh-sec}

Now let $(\CC,F)$ be as above and assume in addition that $F$ is an autoequivalence,
so that $G=F^{-1}$.
Let us take $\Psi=\SS^{-1}$, the inverse Serre functor given by the kernel $\De_{\CC^{op}}^\vee$, 
and use the canonical morphism
$f:F\circ\SS^{-1}\to\SS^{-1}\circ F$ (see \eqref{H-coh-Serre-eq}). 
As we have seen in Corollary \ref{H-coh-cor}, the map
$(F,f)_*$ gets identified with the natural endomorphism $F_*$ of 
the Hochschild cohomology $HH^*(\CC)\simeq\De_{\CC^{op}}^\vee$.
Thus, Theorem \eqref{LR-thm} gives a formula for the supertrace $\sTr(F_*,HH^*(\CC))$.
We will rewrite this formula in a slightly different way.

Note that since $F'$ is an autoequivalence, we have two non-degenerate pairings 
$\Tr(F')\ot\Tr(G)\simeq \Tr(F)\ot\Tr(G')\to k$:
$\lan\cdot,\cdot\ran_{F,G}$ and $\lan\cdot,\cdot\ran_{F',G'}$ (see Section \ref{pairings-sec}).
Let $A_F:\Tr(G)\to \Tr(G)$ be the automorphism such that
$$\lan x,y\ran_{F,G}=\lan A_F(x),y\ran_{F',G'}).$$

\begin{thm} For a tensor autoequivalence $F$ of $\per(\CC)$ one has
$$\sTr(F_*,HH^*(\CC))=\sTr((\SS^{-1},s)_*, \Tr_{\CC}(G)),$$
where
$$s:\SS^{-1}\circ G\to G\circ\SS^{-1}$$
is induced by the isomorphism \eqref{H-coh-Serre-eq}.
Also, we have
$$(\SS^{-1},s)_*=A_F.$$
\end{thm}

\Pf . By Theorem \eqref{LR-thm}, we only have to check that $(\SS^{-1},s)_*=A_F$.
We are going to deduce this from \eqref{2-tau-eq}.
Note that in our case $B=\De_{\CC^{op}}^\vee$, $\Psi=\SS^{-1}$, 
and as the proof of Theorem \ref{LR-thm} shows, $\wt{\b}$ is the map 
$$\SS^{-1}\to \SS^{-1}\circ G\circ F\rTo{s\circ F} G\circ \SS^{-1}\circ F,$$
or on the level of kernels,
$$\De_{\CC^{op}}^\vee\to (F'\Box G)(\De_{\CC^{op}}^\vee).$$
It is easy to see that this map corresponds under the isomorphism of Lemma \ref{F-op-lem}(i)
to $\ga_{F',G'}$. Hence, by Lemma \ref{tau-lem} we obtain
$$\tau^{B}_{F'\Box G}(\wt{\b})=\tau^{\De_{\CC}}_{F\Box G'}(\ga_{F',G'})=\tau_{F',G'}.$$
Thus, \eqref{2-tau-eq} gives the relation
$$\tau_{F',G'}=(\id\ot (\SS^{-1},s)_*)(\tau_{F,G})$$
which is equivalent to what we needed to check.
\ed

In the case $F=\Id$, the above theorem gives the following formula.

\begin{cor}\label{Hoch-coh-hom-cor} One has the equality
$$\sdim HH^*(\CC)=\sTr((\SS^{-1})_*, HH_*(\CC))$$
in $k$.
\end{cor}


\begin{ex} In the situation when $\CC$ is the dg-version of the derived category of coherent sheaves on 
a smooth complex projective variety $X$ the formula of Corollary \ref{Hoch-coh-hom-cor} becomes
$$\sdim HH^*(X)=(-1)^n \sdim HH_*(X),$$
where $n=\dim X$.
Indeed, due to the form of the Serre functor in this case, the operator $(-1)^n\SS^{-1}_*$ on 
$HH_*(X)$ is upper triangular, so its supertrace is equal to $\sdim HH_*(X)$.
Note that the above equality follows also from the well-known identity
$$\sdim HH^*(X)=c_n(\Om^1_X)=(-1)^n c_n(T_X)=(-1)^n \sdim HH_*(X).$$
\end{ex}

\subsection{Example: matrix factorizations}\label{mf-sec}

Let us consider the case when
$\CC=\MF(\w)$ is the category of matrix factorizations of an isolated singularity $\w\in \mg\sub 
R=k[[x_1,\ldots,x_n]]$,
where $k$ is a field of characteristic zero ($\mg$ is the maximal ideal in $R$).
Recall that a matrix factorization $\bar{E}=(E,\de_E)$ of $\w$ is a $\Z/2$-graded
finitely generated free $R$-module $E$ equipped with an odd endomorphism $\de$ such that
$\de^2=\w\cdot\id$. The complexes of morphisms $\Homb(\bar{E}_1,\bar{E}_2)$ are
defined by taking the spaces of $\Z/2$-homogeneous $R$-morphisms $\Hom^*_R(E_1,E_2)$
with the natural differential induced by $\de_{E_1}$ and $\de_{E_2}$.
Note that $\MF(\w)$ is a $\Z/2$-dg-category, however, the formalism developed above works in
a $\Z/2$-graded case as well (cf. \cite[Sec.\ 5.1]{Dyck}).
We will freely use the notation and the conventions of \cite[Sec.\ 2]{PV-HRR}.
Recall that the Hochschild homology of $\w$ is naturally isomorphic to the $\Z/2$-graded space
\begin{equation}\label{H-w-space}
HH_*(\MF(\w))\simeq H(\w):=\Om^n_{R/k}/d\w\we\Om^{n-1}_{R/k}[n]\simeq
(R/(\pa_1\w,\ldots,\pa_n\w))\ot dx_1\we\ldots\we dx_n[n],
\end{equation}
where we set $\pa_i=\pa/\pa x_i$ (see \cite[Thm.\ 6.6]{Dyck}, \cite[(2.28)]{PV-HRR}).

Let 
$$t: (x_1,\ldots,x_n)\mapsto (t_1x_1,\ldots,t_nx_n)$$
be a diagonal symmetry of $\w$, so that
$\w(t_1x_1,\ldots,t_nx_n)=\w(x_1,\ldots,x_n)$.
We associate with $t$ an autoequivalence $t^*$ of $\MF(\w)$ induced by the similar functor
$t^*$ on free $k[[x_1,\ldots,x_n]]$-modules (which acts as identity on objects and as the substitution
$f\mapsto f\circ t$ on morphisms). 

Recall that the kernels giving dg-functors $\MF(\w)\to\MF(\w)$ can be viewed as
matrix factorizations of $\wt{\w}=\w(y_1,\ldots,y_n)-\w(x_1,\ldots,x_n)$ (see \cite[Sec.\ 6.1]{Dyck},
\cite[Sec.\ 4.2]{PV-HRR}).
The diagonal kernel is represented by the stabilized diagonal matrix factorization,
which is the Koszul matrix factorization 
$$\De^{st}=\De^{st}_\w=\{\De_1\w,\ldots,\De_n\w;y_1-x_1,\ldots,y_n-x_n\}=
(R^e\ot {\bigwedge}^\bullet(V),\de^{st})\in \MF(\wt{\w})$$
associated with the decomposition
$$\wt{\w}=\sum_{i=1}^n \De_i\w(x,y)\cdot (y_i-x_i),$$
where $R^e=R\hat{\ot} R=k[[x_1,\ldots,x_n,y_1,\ldots,y_n]]$ and $V=\bigoplus_{i=1}^n k\cdot e_i$
(see \cite[(2.24),(2.25)]{PV-HRR}). The endomorphism $\de^{st}$ of the exterior algebra
has two components: one is given by the exterior product with $\sum \De_i\w\cdot e_i$, and another---by the
contraction with $\sum (y_i-x_i)\cdot e_i^*$.

To write explicitly the kernel corresponding to $t^*$ we use the relation
$$F=(\Id\Box F)(\De_{\CC})$$
for $F=t^*$. The functor $(\Id\Box t^*)$ can be identified with $(\id\times t)^*$, so 
the kernel giving $t^*$ is the Koszul matrix factorization
$(\id\times t)^*\De^{st}$.
Since the functor $\Tr:\MF(\wt{\w})\to \per(k)$ is given by the restriction to the diagonal $y=x$,
we get an isomorphism
$$\Tr(t^*)\simeq \left((\id\times t)^*\De^{st}\right)|_{\De}\simeq \De^{st}|_{\Ga_t},$$
where $\Ga_t$ is the graph of $t$ (given by the equations $y_i=t_ix_i$).
Let us renumber the coordinates in such a way that 
$$t_i\neq 1\ \text{  for } i\le k \ \text{ and }\ t_{k+1}=\ldots=t_n=1.$$
Then the proof of \cite[Lem.\ 2.5.3]{PV-HRR} gives an isomorphism 
\begin{equation}\label{Tr-t*-eq}
H^*\Tr(t^*)\simeq H^*(\De^{st}_{\Ga_t})\simeq H(\w_t),
\end{equation}
where $H(\w_t)$ is the space \eqref{H-w-space} for the potential
$\w_t$, the restriction of $\w$ to the subspace $x_1=\ldots=x_k=0$ of fixed points of $\w$.
Note that $\w_t$ still has isolated singularity (see \cite[Lem.\ 2.5.3]{PV-HRR}(i)).
More explicitly, the isomorphism \eqref{Tr-t*-eq} is induced by the projection from the
term $R^e\ot\bigwedge^{n-k}(V)$ to $R/(x_1,\ldots,x_k)\ot e_{k+1}\we\ldots\we e_n$ 
(which factors through the restriction to $\De\cap\Ga_t$).

Now let $\bar{E}=(E,\de_E)$ be a matrix factorization of $\w$.
The computation of the generalized boundary-bulk map
\begin{equation}\label{mf-bb-map}
\tau^{\bar{E}}_{t^*}:\Hom(\bar{E},t^*\bar{E})\to\Tr(t^*)\simeq H(\w_t)
\end{equation}
is very similar to the computation of the equivariant Chern character in \cite[Sec.\ 3]{PV-HRR}.
By definition, we have to take the restriction to the diagonal of the morphism 
$c^{\bar{E}}_{t^*}=(\id\times t)^*c^{\bar{E}}$, where
$$c^{\bar{E}}=c^{\bar{E}}_{\Id}:\bar{E}^*\boxtimes\bar{E}\to\De^{st}$$
is the canonical morphism \eqref{A-vee-Id-eq}. The latter morphism was explicitly described in
\cite[Sec.\ 3.1, 3.2]{PV-HRR}. It is convenient to use the identification of the complex of morphisms in 
$\MF(\w)$ from $\bar{E}^*\boxtimes\bar{E}$ to $\De^{st}$ with the tensor product
$$L:=\De^{st}\ot_{R^e}(\bar{E}\boxtimes\bar{E}^*)$$
(see \cite[(3.1)]{PV-HRR}). Then the morphism $c^{\bar{E}}$ corresponds to a certain
even closed element $D\in L$. Let us choose a trivialization
$E\simeq U\ot R$ of the $R$-module $E$, where $U$ is a $\Z/2$-graded vector space.
This leads to the identification 
$$L\simeq \bigoplus_{j=0}^n {\bigwedge}^j(V)\ot \End(U)\ot R^e$$
and writing 
$$D=\sum_j D_j,\ \text{ with }\ D_j=\sum_{i_1<\ldots<i_j}e_{i_1}\we\ldots\we e_{i_j}\ot D_j(i_1,\ldots,i_j)$$
with the components 
$D_j(i_1,\ldots,i_j)\in \End(U)\ot R^e$ and $D_0=1\ot\id_U$, we obtain a system of equations on 
these components with the unique solution such that $D_j(i_1,\ldots,i_j)$ does not depend on $y_r$
with $r<i_1$ (see \cite[Lem.\ 3.2.1]{PV-HRR}). 
The morphism $(\id\times t)^*c^{\bar{E}}$ corresponds to the element
$$(\id\times t)^*D\in (\id\times t)^*L\simeq {\bigwedge}^\bullet(V)\ot \End(U)\ot R^e,$$
where the latter identification uses the action of $\id\times t$ on $R^e$.

As the description of the isomorphism \eqref{Tr-t*-eq} shows, 
the map $\tau^{\bar{E}}_{t^*}=\Tr((\id\times t)^*c^{\bar{E}})$ is determined by 
$$(\id\times t)^*D_{n-k}(k+1,\ldots,n)|_{\De\cap\Ga_t}=
\pa_n\de_E(x)\circ\ldots\circ\pa_{k+1}\de_E(x)|_{x_1=\ldots=x_k=0},$$
where we view $\de_E$ as an element of $\End(U)\ot R$ 
(see the proof of \cite[Thm.\ 3.3.3]{PV-HRR}).
Thus, if use the identification $\Homb(\bar{E},t^*\bar{E})\simeq (\End(U)\ot R,d)$, 
where $d$ is the differential induced by $\de_E$, then we arrive to the formula for
$\tau^{\bar{E}}_{t^*}$ similar to those of \cite[Thm.\ 3.2.3, Thm.\ 3.3.3]{PV-HRR}.

\begin{prop}\label{mf-gen-bb-prop} 
In the above situation one has
$$\tau^{\bar{E}}_{t^*}(\a)=\sTr\left([\pa_n\de_E\circ\ldots\circ \pa_{k+1}\de_E\circ\a]|_{x_1=\ldots=x_k=0}\right)
\cdot dx_{k+1}\we\ldots\we dx_n \ \mod\ (\pa_{k+1}\w_t,\ldots,\pa_n\w_t).$$
\end{prop}

Next, we want to compute the canonical pairing \eqref{Tr-F-pairing} between $\Tr(F)$ and $\Tr(G)$ for 
the autoequivalence $F=t^*$ of $\MF(\w)$. Note that in this case $G=(t^{-1})^*$, so
both $\Tr(F)$ and $\Tr(G)$ are naturally isomorphic to the $\Z/2$-graded space $H(\w_t)=H(\w_{t^{-1}})$.
Again the calculation is very similar to that of the canonical pairing on the equivariant Hochschild
homology in \cite[Sec.\ 4.2]{PV-HRR}. Namely, it is enough to compute the corresponding Casimir
tensor $\tau_{F,G}\in \Tr(F')\ot\Tr(G)$, which is obtained by applying the generalized boundary-bulk
map associated with $F'\Box G$ to the canonical morphism $\De\to (F'\Box G)(\De)$ 
(see Lemma \ref{Casimir-lem}(ii)).
In our case $F'$ is the autoequivalence $(t^{-1})^*$ of $\MF(-\w)$, so 
$F'\Box G$ is the autoequivalence $(t^{-1}\times t^{-1})^*$ of $\MF(\wt{\w})$.
Thus, to calculate $\tau_{F,G}$ we can apply Proposition \ref{mf-gen-bb-prop} to the canonical
isomorphism
\begin{equation}\label{a-De-st-isom}
\a:\De^{st}\to (t^{-1}\times t^{-1})^*\De^{st}
\end{equation}
in the homotopy category of $\MF(\wt{\w})$. Note that
$$(t^{-1}\times t^{-1})^*\De^{st}\simeq \{(t^{-1}\times t^{-1})^*\De_1\w,\ldots,(t^{-1}\times t^{-1})^*\De_n\w;
t_1^{-1}(y_1-x_1),\ldots,t_n^{-1}(y_n-x_n)\}$$
and the morphism $\a$ is characterized uniquely (up to homotopy)
by the condition that the induced map on $\coker(\de_{10}:\bigwedge^1\to\bigwedge^0)$
is identity, where $\de_{10}$ is the component of the differential $\de$ of the Koszul matrix factorization
(see \cite[Prop.\ 2.3.1]{PV-HRR} or \cite{Dyck}). We will use the following fact about Koszul matrix factorizations
(that follows from the proof of \cite[Lem.\ 2.5.5]{PV-HRR}).

\begin{lem} Let $b_\bullet=(b_1,\ldots,b_r)$ be a regular sequence in a local ring $A$,
and let
$$a_1b_1+\ldots+a_rb_r=a'_1b_1+\ldots+a'_rb_r=\w\in A$$
for some sequences $a_\bullet$ and $a'_\bullet$ in $A$. Then
there exists an isomorphism of matrix factorizations
$$\{a_\bullet,b_\bullet\}\rTo{\sim}\{a'_\bullet,b_\bullet\}$$
of the form $\exp(h)\we ?$ for some $h\in\bigwedge^2_A(A^r)$.
\end{lem}

Let us denote by $\a_t$ the automorphism of the exterior algebra $\bigwedge^\bullet(V)$ induced by
the map $t:V\to V$ sending $e_i$ to $t_i\cdot e_i$. Then 
$\a_t$ gives an isomorphism of Koszul matrix factorizations
$$\a_t: \{t_1^{-1}f_1,\ldots,t_n^{-1}f_n;(y_1-x_1),\ldots,(y_n-x_n)\}\to
\{f_1,\ldots,f_n;t_1^{-1}(y_1-x_1),\ldots,t_n^{-1}(y_n-x_n)\},$$
where we denoted $f_i(x,y)=(t^{-1}\times t^{-1})^*\De_i\w$.
Now applying the above lemma we get an isomorphism \eqref{a-De-st-isom} of the form
$\a=\a_t\circ \exp(h)$. Note that $\a$ preserves the filtration
$$F_p={\bigwedge}^{\ge p}(\bigoplus_{i=1}^r k\cdot e_i)\ot
{\bigwedge}^{\bullet}(\bigoplus_{j=r+1}^n k\cdot e_j)\ot R^e$$
on ${\bigwedge}^{\bullet}(V)\ot R^e$ (since $F_p$ are ideals with respect to the exterior multiplication). 
Thus, the same argument as in \cite[Sec.\ 4.2]{PV-HRR} shows the equality
$$\tau_{F,G}=\ch(\De^{st}_{\w_t})\cdot\det[\id-t,V/V^t]=\ch(\De^{st}_{\w_t})\cdot (1-t_1)\ldots (1-t_k),$$
where $\De^{st}_{\w_t}$ is the stabilized diagonal for the potential $\w_t(x_{k+1},\ldots,x_n)$.
Hence, using \cite[Prop.\ 4.1.1, 4.1.2]{PV-HRR} we arrive to the following result.

\begin{prop}\label{mf-pairing-prop}
Under the isomorphisms $\Tr(t^*)\simeq\Tr((t^{-1})^*)\simeq H(\w_t)$ the canonical pairing
$\lan\cdot,\cdot\ran_{t^*,(t^{-1})^*}$ (see \eqref{Tr-F-pairing}) is given by
$$\lan\cdot,\cdot\ran_{t^*,(t^{-1})^*}=(1-t_1)^{-1}\ldots(1-t_k)^{-1}\cdot\lan\cdot,\cdot\ran_{\w_t},$$
where $\lan\cdot,\cdot\ran_{\w_t}$ is the canonical pairing on $H(\w_t)=HH_*(\MF(\w_t))$, given by
\cite[(4.6)]{PV-HRR}.
\end{prop}

In the case when $0$ is an isolated fixed point of $t$, by Propositions
\ref{mf-gen-bb-prop} and \ref{mf-pairing-prop}, the holomorphic Lefschetz formula \eqref{HLF-eq}
takes the form stated in Theorem \ref{mf-isolated-thm}.

\begin{rem} Comparing the formula of Proposition \ref{mf-gen-bb-prop} 
with that of \cite[Thm.\ 3.3.3]{PV-HRR} we see
that in the situation when $G$ is a finite group of symmetries of $\w$, and $\bar{E}$ is a $G$-equivariant
matrix factorization of $\w$, the components
of the boundary-bulk map $\tau^{\bar{E}}$ 
for the dg-category $\MF_G(\w)$ of $G$-equivariant matrix factorizations of $\w$ are given by
$$\a\mapsto\tau^{\bar{E}}_{g^*}(g\circ\a), \ \ g\in G,$$
where we use the action of $g\in G$ on $E$. Similarly, comparing Proposition \ref{mf-pairing-prop}
with the formula for the canonical pairing on the Hochschild homology of $\MF_G(\w)$, we see that
the latter pairing has components 
$$|G|^{-1}\lan\cdot,\cdot\ran_{g^*,(g^{-1})^*}.$$
Elsewhere we will show that the same relations hold whenever a finite group $G$ acts on a smooth
and compact dg-category and we pass to the corresponding $G$-equivariant category. 
\end{rem}

In conclusion let us point out several consequences of
the Lefschetz formula \eqref{HLF-eq} in the situation considered above.

\begin{cor}\label{zero-cor} 
Let $\w(x_1,\ldots,x_n)$ be an isolated singularity, $t\in (k^*)^n$ a diagonal symmetry of
$\w$, $\bar{A}$ and $\bar{B}$ matrix factorizations of $\w$ equipped with closed morphisms
of degree zero $\a:\bar{A}\to t^*\bar{A}$ and $\b:t^*\bar{B}\to\bar{B}$.
Assume that the fixed locus of $t$ has odd dimension (i.e., the number of coordinates of $t$
that are equal to $1$ is odd). Then the endomorphism 
$(t^*,\a,\b)_*$ of $\Hom(\bar{A},\bar{B})$ defined by \eqref{F-a-b-eq} has zero supertrace.
\end{cor}

\Pf . The map $\tau^{\bar{A}}_{t_*}$ takes values in the $\Z/2$-graded space 
$\Tr(t^*)\simeq H(\w_t)$ which lives in odd degree due to our assumption on the fixed locus of $t$.
Hence, $\tau^{\bar{A}}_{t_*}(\a)=0$ and the right-hand side of the formula \eqref{HLF-eq}
vanishes in our situation.
\ed

Note that Corollary \ref{zero-cor} generalizes the result conjectured by Dao \cite{Dao} and proved
in \cite{PV-HRR} and \cite{B-vS} (see also \cite{MPSW}) that 
$\Hom(\bar{A},\bar{B})$ has zero Euler characteristic in the case when $n$ is odd.

\begin{prop}\label{divisibility-cor} Assume $\operatorname{char}(k)=0$.
Suppose we have an isolated singularity $\w(x_1,\ldots,x_n)\in k[[x_1,\ldots,x_n]]$ with a diagonal symmetry
$t=(t_1,\ldots,t_n)\in (k^*)^n$ of prime order $p$, such that
$t_i\neq 1$ for all $i$ (i.e., the origin is an isolated fixed point of $t$). 
Let $\bar{A}$ be a $\Z/p$-equivariant matrix factorization of $\w$. Then the complex
of $\Z/p$-modules $\bar{A}|_0$ defines a class $[\bar{A}|_0]$ in the representation ring $R_{\Z/p}$. 
Assume that $\left \lceil \frac{n}{2} \right  \rceil -1\ge m(p-1)$. Then
$[\bar{A}|_0]$ is divisible by $p^m$ in $R_{\Z/p}$.
\end{prop}

\Pf . The formula \eqref{HLF-eq} applied to $\bar{B}=\bar{A}$ gives an equality 
\begin{equation}\label{trace-A-A-eq}
\sTr(\a|_0,\bar{A}|_0)\cdot\sTr(\a^{-1}|_0,\bar{A}|_0)=\sTr(t^*,\Hom(\bar{A},\bar{A}))
\cdot\prod_{i=1}^n (1-t_i)
\end{equation}
in $k$, where $\a:\bar{A}\rTo t^*\bar{A}$ is the $\Z/p$-equivariant structure on $\bar{A}$.
We can view this as an equality in $\Z[\zeta_p]\sub\ov{k}$, where $\zeta_p$ is the primitive $p$th root of unity.
Let $\tr:R_{\Z/p}\to\Z[\zeta_p]$ be the homomorphism induced by taking the trace of the generator of $\Z/p$.
We claim that it restricts to an isomorphism of the ideal $I\sub R_{\Z/p}$ of virtual representations of
rank $0$ with the ideal $(1-\zeta_p)\sub \Z[\zeta_p]$. Indeed, a class
$\sum_{i=0}n_i\chi^i\in R_{\Z/p}$ belongs to $I$ if $\sum_{i=0}^{p-1}n_i=0$. In this case 
$$\sum_{i=0}^{p-1}n_i\zeta_p^i=(\zeta_p-1)\cdot\sum_{i=1}^{p-1}n_i\frac{\zeta_p^i-1}{\zeta_p-1},$$
and our claim follows from the fact that $(\frac{\zeta_p^i-1}{\zeta_p-1})_{1\le i\le p-1}$ is
a basis of $\Z[\zeta_p]$ over $\Z$. Note that the class $[\bar{A}|_0]\in R_{\Z/p}$ belongs to $I$ (since $\w$
is not a zero divisor). 
Let $a=\tr([\bar{A}|_0])\in (1-\zeta)\Z[\zeta_p]$. It is enough to prove that under our assumptions
$a/(1-\zeta_p)$ is divisible by $p^m$ in $\Z[\zeta_p]$.
The formula \eqref{trace-A-A-eq} implies that
$a\cdot\ov{a}$ is divisible by $(1-\zeta_p)^n$ in $\Z[\zeta_p]$. 
Since $\Nm_{\Z[\zeta_p]/\Z}(1-\zeta_p)=p$, the ideal
$(1-\zeta_p)\sub \Z[\zeta_p]$ is prime. This ideal is also stable under conjugation, so we derive
that $a$ is divisible by $(1-\zeta_p)^{\left \lceil \frac{n}{2} \right  \rceil}$. Hence, 
$a/(1-\zeta_p)$ is divisible by $(1-\zeta_p)^{\left \lceil \frac{n}{2} \right  \rceil-1}$
It remains to use the fact that
$$p=\prod_{i=1}^{p-1}(1-\zeta_p^i)=u\cdot (1-\zeta_p)^{p-1}$$
where $u$ is a unit in $\Z[\zeta_p]$, so we can replace $(1-\zeta_p)^{m(p-1)}$ by $p^m$.
\ed

\begin{ex}\label{divisibility-ex}
Assume that an isolated singularity $\w(x_1,\ldots,x_n)$ is even, i.e.,
$\w(-x)=\w(x)$. Then Proposition \ref{divisibility-cor} states (in the case $p=2$) that for any $\Z/2$-equivariant
matrix factorization $\bar{A}$ of $\w$, the class $[\bar{A}|_0]$ is divisible by
$2^{\left \lceil \frac{n}{2} \right  \rceil -1}$ in $R_{\Z/2}$. Equivalently, we have the divisibility
\begin{equation}\label{main-divisibility-eq}
2^{\left \lceil \frac{n}{2} \right  \rceil}|\sTr((-1)^*,\bar{A}|_0)
\end{equation}
(see the proof of Proposition \ref{divisibility-cor}).
Thus, either $[\bar{A}|_0]=0$ in $R_{\Z/2}$ or
\begin{equation}\label{even-inequality}
\rk(\bar{A})=\dim_k(\bar{A}|_0)\ge 2^{\left \lceil \frac{n}{2} \right  \rceil}.
\end{equation}
\end{ex}

\begin{rems}
1. Recall that one conjectures for a nonzero finitely generated module $M$ of finite projective dimension over a 
local noetherian commutative ring $R$ the following lower bound on the total Betti number of $M$:
$$\sum_i \dim_k \Tor^R_i(M,k)\ge 2^{d(R)-d(M)},$$
where $d(\cdot)$ is the Krull dimension. 
In \cite{AB-bounds} this inequality is proved in the case when $R$ is a graded commutative algebra $R$, finitely generated by $R_1$ over $R_0=k$,
under certain additional assumption on a graded $R$-module $M$ (e.g., if $M$ has odd multiplicity).
Our proof of Proposition \ref{divisibility-cor} is very similar to that of Theorem 3 of \cite{AB-bounds} 
which essentially proves certain divisibility using Hilbert series and then deduces the lower bound from it.
 
2. Set $R=k[[x_1,\ldots,x_n]]$. Recall that there is an equivalence between the stable category of $S=R/(\w)$
and the homotopy category of matrix factorizations of $\w$ (see \cite{Bu, Orlov}), and that
the matrix factorization $\bar{A}=(A,\de)$ corresponding to a finitely generated $S$-module $M$ 
specializes to the $2$-periodic free $S$-resolution
of a maximal Cohen-Macaulay module stably equivalent to $M$ (cf. \cite{Eis, Orlov}). In particular,
we get
$$H^*(A|_0,\de)\simeq\Tor_{2N}^S(M,k)\oplus\Tor_{2N+1}^S(M,k).$$ 
Thus, from \cite[Conj.\ 7.5]{AB-rs} we would get the following bound on the rank
of the matrix factorization $\bar{A}$ of $\w(x_1,\ldots,x_n)$ corresponding
to a nonzero $S$-module $M$ over $S$ of infinite projective dimension:
$$\rk(\bar{A})\ge 2^{n-d(M)}$$
(to deduce this from \cite[Conj.\ 7.5]{AB-rs} one has to consider an appropriate regular sequence
annihilating $M$). Note that 
\eqref{even-inequality} gives a stronger bound on $\rk(\bar{A})$ provided
$d(M)>\left \lfloor \frac{n}{2} \right \rfloor$ (but this bound needs additional assumptions).
In the graded case this can be interpreted in terms of the multiplicity polynomial---see Corollary below.
\end{rems}

\begin{cor}\label{2nd-divisibility-cor}  Assume $\operatorname{char}(k)=0$.
Let $M$ be a graded module over $R=k[x_1,\ldots,x_n]$ (where $\deg(x_i)=1$).
Let $e_M(t)$ be the multiplicity polynomial of $M$, so that the Hilbert series of $M$ has
form $H_M(t)=e_M(t)(1-t)^{-d(M)}$, where $d(M)$ is the Krull dimension of $M$.
Assume that $M$ is annihilated by a homogeneous polynomial of even degree $\w$ such that
the corresponding projective hypersurface is smooth. Then
$e_M(-1)$ is divisible by $2^{d(M)-\left \lfloor \frac{n}{2} \right \rfloor}$.
\end{cor}

\Pf . Let us consider the polynomial
$$\chi_M(t)=\sum_i (-1)^i H_{\Tor^R_i(M,k)}(t).$$
Since $H_R(t)=(1-t)^{-n}$, we get
$$\chi_M(t)=H_M(t)(1-t)^n,$$
or in terms of the multiplicity polynomial $e_M(t)$,
\begin{equation}\label{chi-e-M-eq}
\chi_M(t)=e_M(t) (1-t)^{n-d(M)}.
\end{equation}
On the other hand, we can view $M$ as a $\Z/2$-equivariant module over $R/(\w)$.
From the graded free resolution of $M$ over $R$ we get by stabilization a $\Z/2$-equivariant
matrix factorization $\bar{A}$ of $\w$ with 
$$\sTr((-1)^*,\bar{A}|_0)=\chi_M(-1).$$
Now the required divisibility for $e_M(-1)$ is obtained by substituting $t=-1$ into \eqref{chi-e-M-eq} 
and using \eqref{main-divisibility-eq}.
\ed

\end{document}